# SHARP LIPSCHITZ ESTIMATES FOR OPERATOR $\bar\partial_{\mathbf{M}}$ ON A $q$-PSEUDOCONCAVE CR MANIFOLD

PETER L. POLYAKOV

ABSTRACT. We construct integral operators $R_r$ and $H_r$ on a regular q-pseudoconcave CR manifold $\mathbf{M}$ such that
$$f = \bar\partial_{\mathbf{M}} R_r(f) + R_{r+1}(\bar\partial_{\mathbf{M}} f) + H_r(f),$$
for $f \in C^\infty_{(0,r)}(\mathbf{M})$ and prove sharp estimates in a special Lipschitz scale.

## 1. INTRODUCTION.

Let $\mathbf{M}$ be a CR submanifold in a complex $n$ - dimensional manifold $\mathbf{G}$ such that for any $z \in \mathbf{M}$ there exist a neighborhood $V \ni z$ in $\mathbf{G}$ and smooth real valued functions
$$\{\rho_k, \ k=1,\ldots,m \ (1 < m < n-1)\}$$
on $V$ such that
$$\begin{aligned} \mathbf{M} \cap V &= \{z \in \mathbf{G} \cap V : \rho_1(z) = \cdots = \rho_m(z) = 0\}, \\ \partial\rho_1 &\wedge \cdots \wedge \partial\rho_m \neq 0 \ \text{ on } \mathbf{M} \cap V. \end{aligned} \quad (1)$$

In this paper we continue the study of regularity of the operator $\bar\partial_{\mathbf{M}}$ on a submanifold $\mathbf{M}$ satisfying special concavity condition. In [P] we considered solutions of the $\bar\partial_{\mathbf{M}}$ equation with an $L^\infty$ right hand side. Here we prove sharp estimates in a special Lipschitz scale.

Before formulating the main result we will introduce necessary notations and definitions.

The CR structure on $\mathbf{M}$ is induced from $\mathbf{G}$ and is defined by the subbundles
$$T''(\mathbf{M}) = T''(\mathbf{G})|_{\mathbf{M}} \cap \mathbf{C}T(\mathbf{M}) \quad \text{and} \quad T'(\mathbf{M}) = T'(\mathbf{G})|_{\mathbf{M}} \cap \mathbf{C}T(\mathbf{M}),$$
where $\mathbf{C}T(\mathbf{M})$ is the complexified tangent bundle of $\mathbf{M}$ and the subbundles $T''(\mathbf{G})$ and $T'(\mathbf{G}) = \overline{T''(\mathbf{G})}$ of the complexified tangent bundle $\mathbf{C}T(\mathbf{G})$ define the complex structure on $\mathbf{G}$.

We will denote by $T^c(\mathbf{M})$ the subbundle $T(\mathbf{M}) \cap [T'(\mathbf{M}) \oplus T''(\mathbf{M})]$. If we fix a hermitian scalar product on $\mathbf{G}$ then we can choose a subbundle $N \in T(\mathbf{M})$ of real dimension $m$ such that $T^c(\mathbf{M}) \perp N$ and for a complex subbundle $\mathbf{N} = \mathbf{C}N$ of $\mathbf{C}T(\mathbf{M})$ we have
$$\mathbf{C}T(\mathbf{M}) = T'(\mathbf{M}) \oplus T''(\mathbf{M}) \oplus \mathbf{N}, \ \ T'(\mathbf{M}) \perp \mathbf{N} \ \text{ and } \ T''(\mathbf{M}) \perp \mathbf{N}.$$

We define the Levi form of $\mathbf{M}$ as the hermitian form on $T'(\mathbf{M})$ with values in $\mathbf{N}$
$$\mathcal{L}_z(L(z)) = \sqrt{-1} \cdot \pi\left(\left[\overline{L}, L\right]\right)(z) \quad (L(z) \in T'_z(\mathbf{M})),$$
where $\left[\overline{L}, L\right] = \overline{L}L - L\overline{L}$ and $\pi$ is the projection of $\mathbf{C}T(\mathbf{M})$ along $T'(\mathbf{M}) \oplus T''(\mathbf{M})$ onto $N$.

If functions $\{\rho_k\}$ are chosen so that the vectors $\{\text{grad}\rho_k\}$ are orthonormal then the Levi form of $\mathbf{M}$ may be defined as
$$L_z(\mathbf{M}) = -\sum_{k=1}^m (L_z\rho_k(\zeta)) \cdot \text{grad } \rho_k(z),$$

Date: April 28, 1998.
1991 *Mathematics Subject Classification.* 32F20, 32F10.
The author was partially supported by NSF Grant DMS-9022140 during his visit to MSRI.

where $L_z\rho(\zeta)$ is the Levi form of the real valued function $\rho \in C^4(\mathbf{D})$ at the point $z$:

$$L_z\rho(\zeta) = \sum_{i,j} \frac{\partial^2 \rho}{\partial \zeta_i \partial \bar\zeta_j}(z)\, \zeta_i \cdot \bar\zeta_j.$$

For a pair of vectors $\mu = (\mu_1, \ldots, \mu_n)$ and $\nu = (\nu_1, \ldots, \nu_n)$ in $\mathbb{C}^n$ we will denote $\langle \mu, \nu \rangle = \sum_{i=1}^n \mu_i \cdot \nu_i$.

For a unit vector $\theta = (\theta_1, \ldots \theta_m) \in \mathrm{Re}\mathbf{N}_z$ we define the Levi form of $\mathbf{M}$ at the point $z \in \mathbf{M}$ in the direction $\theta$ as the scalar hermitian form on $\mathbf{C}T_z^c(\mathbf{M})$

$$\langle \theta,\ L_z(\mathbf{M}) \rangle = -L_z\rho_\theta(\zeta),$$

where $\rho_\theta(\zeta) = \sum_{k=1}^m \theta_k \rho_k(\zeta)$.

Following [H2] we introduce the notion of a q-pseudoconcave CR manifold. Namely, we call $\mathbf{M}$ q-pseudoconcave (weakly q-pseudoconcave) at $z \in \mathbf{M}$ in the direction $\theta$ if the Levi form of $\mathbf{M}$ at $z$ in this direction $\langle \theta,\ L_z(\mathbf{M}) \rangle$ has at least $q$ negative ($q$ nonpositive) eigenvalues on $\mathbf{C}T_z^c(\mathbf{M})$.

We call $\mathbf{M}$ q-pseudoconcave (weakly q-pseudoconcave) at $z \in \mathbf{M}$ if it is q-pseudoconcave (weakly q-pseudoconcave) in all directions.

We call a q-pseudoconcave CR manifold $\mathbf{M}$ by a regular q-pseudoconcave CR manifold (cf. [P]) if for any $z \in \mathbf{M}$ there exist an open neighborhood $\mathcal{U} \ni z$ in $\mathbf{M}$ and a family $E_q(\theta, z)$ of q-dimensional complex linear subspaces in $\mathbf{C}T_z^c(\mathbf{M})$ smoothly depending on $(\theta, z) \in \mathbb{S}^{m-1} \times \mathcal{U}$ and such that the Levi form $\langle \theta,\ L_z(\mathbf{M}) \rangle$ is strictly negative on $E_q(\theta, z)$.

Following [S] we define spaces $\Gamma^\beta(\mathbf{M})$ for $0 < \beta < 2$ with the norm

$$\|h\|_{\Gamma^\beta(\mathbf{M})} = \|h\|_{\Lambda^{\frac{\beta}{2}}(\mathbf{M})} + \sup\left\{\|h(x(\cdot))\|_{\Lambda^\beta([0,1])}\right\}$$

where the $\sup$ is taken over all curves $x: [0,1] \to \mathbf{M}$ such that

$$\begin{aligned} |x'(s)|, |x''(s)| &\leq 1, \\ x'(s) &\in T^c(\mathbf{M}). \end{aligned} \tag{2}$$

We introduce spaces $\Pi^a(\mathbf{M})$ for positive $a = p + \alpha$ with $p \in \mathbb{Z}$ and $0 < \alpha < 1$ by saying that function $h \in \Pi^a(\mathbf{M})$ if

(a) for any set of tangent vector fields $D_1, \ldots, D_k, D_1^c, \ldots, D_s^c$ on $\mathbf{M}$ such that

$$\|D_i\|_{C^{p+1}(\mathbf{M})}, \|D_i^c\|_{C^{p+1}(\mathbf{M})} \leq 1$$

with $D_i^c \in \mathbf{C}T^c(\mathbf{M})$ and $2k + s \leq p$

$$\|D_1^c \circ \cdots \circ D_s^c \circ D_1 \circ \cdots \circ D_k h\|_{\Gamma^\alpha(\mathbf{M})} < \infty,$$

(b) for any set of tangent vector fields $D_1, \ldots, D_k, D_1^c, \ldots, D_s^c$ on $\mathbf{M}$ such that

$$\|D_i\|_{C^p(\mathbf{M})}, \|D_i^c\|_{C^p(\mathbf{M})} \leq 1$$

with $D_i^c \in \mathbf{C}T^c(\mathbf{M})$ and $2k + s \leq p - 1$

$$\|D_1^c \circ \cdots \circ D_s^c \circ D_1 \circ \cdots \circ D_k h\|_{\Gamma^{1+\alpha}(\mathbf{M})} < \infty.$$

Accordingly we define

$$\|h\|_{\Pi^a(\mathbf{M})} = \sup_{2k+s \leq p} \|\overbrace{D^c \circ}^{s} \overbrace{D}^{k} h\|_{\Gamma^\alpha(\mathbf{M})} + \sup_{2k+s \leq p-1} \|\overbrace{D^c \circ}^{s} \overbrace{D}^{k} h\|_{\Gamma^{1+\alpha}(\mathbf{M})}.$$

For a differential form $g = \sum_{I,J} g_{I,J}(z) dz^I \wedge d\bar z^J$ with $|I| = k$ and $|J| = r$ we say that $g \in \Pi_{(k,r)}^a(\mathbf{M})$ if $g_{I,J} \in \Pi^a(\mathbf{M})$.

The following theorem represents the main result of the paper.

**Theorem 1.** *Let $a > 1$ and let a compact $C^\infty$ submanifold $\mathbf{M} \subset \mathbf{G}$ of the form (1) be regular $q$-pseudoconcave. Then for any $r = 1, \ldots, q-1$ there exist linear operators*

$$\mathbf{R}_r : \Pi^a_{(0,r)}(\mathbf{M}) \to \Pi^{a+1}_{(0,r-1)}(\mathbf{M}) \quad and \quad \mathbf{H}_r : \Pi^a_{(0,r)}(\mathbf{M}) \to \Pi^a_{(0,r)}(\mathbf{M})$$

*such that $\mathbf{R}_r$ is bounded and $\mathbf{H}_r$ is compact and such that for any differential form $f \in C^\infty_{(0,r)}(\mathbf{M})$ the equality:*

$$f = \bar{\partial}_{\mathbf{M}} \mathbf{R}_r(f) + \mathbf{R}_{r+1}(\bar{\partial}_{\mathbf{M}} f) + \mathbf{H}_r(f) \tag{3}$$

*holds.*

In [P] the existence of operators $R_r : L^\infty_{(0,r)}(\mathbf{M}) \to \Gamma^1_{(0,r-1)}(\mathbf{M})$ and $H_r : L^\infty_{(0,r)}(\mathbf{M}) \to L^\infty_{(0,r-1)}(\mathbf{M})$, satisfying (3) was proved. Barrier functions $\Phi(\zeta, z)$ used in [P] allow to prove compactness of $H_r$ only for $1 \leq r < q - m + 1$, not for $1 \leq r < q$ as it is mistakenly stated in the Proposition 6 of [P].

Here we use different barrier functions, which are closer to the barrier functions in [AiH] and have such a property that corresponding "Cauchy terms" in the integral formulas disappear for $1 \leq r < q$.

A version of the main theorem for q-pseudoconvex hypersurfaces ($m = 1$) and spaces $\Gamma^a(\mathbf{M})$ with $a \geq 0$ was proved by G. B. Folland and E. M. Stein in [FS] (cf. also [H1]).

For $q$-pseudoconcave manifolds of codimension higher than one I. Naruki in [Na] using Kohn-Hormander's method constructed bounded operators $R_r : L^2_{(0,r)}(\mathbf{M}) \to L^2_{(0,r-1)}(\mathbf{M}')$ for the $(0, r)$ forms with $r > n - m - q$.

Then in [H2] and [AiH] with the use of explicit integral formulas bounded operators $R_r : L^\infty_{(0,r)}(\mathbf{M}) \to \Gamma^{0,1-\epsilon}_{(0,r-1)}(\mathbf{M})$ were constructed on a q-pseudoconcave CR manifold of higher codimension for the forms of type $(0, r)$ with $r < q$ or $r > n - m - q$.

Existence of a solution $f$ of equation $\bar{\partial}_{\mathbf{M}} f = g$ such that $D^c f \in \Gamma^\alpha_{(0,r-1)}(\mathbf{M})$ for $D^c \in \mathbf{C}T^c(\mathbf{M})$ and $g \in \Gamma^\alpha_{(0,r)}(\mathbf{M})$, $0 < \alpha < 1$ and $\mathbf{M}$ - quadratic q-pseudoconcave CR manifold was obtained in the paper [BGG] by R. Beals, B. Gaveau and P.C. Greiner.

Author thanks G. Henkin for helpful discussions.

## 2. Construction of $\mathbf{R}_r$ and $\mathbf{H}_r$.

For a vector-valued function $\eta = (\eta_1, \ldots, \eta_n)$ we will use the notation:

$$\omega'(\eta) = \sum_{k=1}^n (-1)^{k-1} \eta_k \wedge_{j \neq k} d\eta_j, \quad \omega(\eta) = \wedge_{j=1}^n d\eta_j.$$

If $\eta = \eta(\zeta, z, t)$ is a smooth function of $\zeta \in \mathbb{C}^n, z \in \mathbb{C}^n$ and a real parameter $t \in \mathbb{R}^p$ satisfying the condition

$$\sum_{k=1}^n \eta_k(\zeta, z, t) \cdot (\zeta_k - z_k) = 1 \tag{4}$$

then

$$d\omega'(\eta) \wedge \omega(\zeta) \wedge \omega(z) = 0$$

or, separating differentials,

$$d_t \omega'(\eta) + \bar{\partial}_\zeta \omega'(\eta) + \bar{\partial}_z \omega'(\eta) = 0. \tag{5}$$

Also, if $\eta(\zeta, z, t)$ satisfies (4) then the differential form $\omega'(\eta) \wedge \omega(\zeta) \wedge \omega(z)$ can be represented as:
$$\sum_{r=0}^{n-1} \omega'_r(\eta) \wedge \omega(\zeta) \wedge \omega(z), \tag{6}$$
where $\omega'_r(\eta)$ is a differential form of the order $r$ in $d\bar{z}$ and respectively of the order $n-r-1$ in $d\bar{\zeta}$ and $dt$. From (5) and (6) follow equalities:
$$d_t \omega'_r(\eta) + \bar{\partial}_\zeta \omega'_r(\eta) + \bar{\partial}_z \omega'_{r-1}(\eta) = 0 \qquad (r = 1, \ldots, n), \tag{7}$$
and
$$\omega'_r(\eta) = \frac{1}{(n-r-1)!r!} \mathrm{Det}\left[\eta,\ \overbrace{\bar{\partial}_z \eta}^{r},\ \overbrace{\bar{\partial}_{\zeta,t} \eta}^{n-r-1}\right], \tag{8}$$
where the determinant is calculated by the usual rules but with external products of elements and the position of the element in the external product is defined by the number of its column.

Let $\widetilde{\mathcal{U}}$ be an open neighborhood in $\mathbf{G}$ and $\mathcal{U} = \widetilde{\mathcal{U}} \cap \mathbf{M}$. We call a vector function
$$P(\zeta, z) = (P_1(\zeta, z), \ldots, P_n(\zeta, z)) \quad \text{for} \quad (\zeta, z) \in \left(\widetilde{\mathcal{U}}\right) \times \widetilde{\mathcal{U}}$$
by strong $\mathbf{M}$-barrier for $\widetilde{\mathcal{U}}$ if there exists $C > 0$ such that the inequality:
$$|\Phi(\zeta, z)| > C \cdot \left(\rho(\zeta) + |\zeta - z|^2\right) \tag{9}$$
holds for $(\zeta, z) \in \left(\widetilde{\mathcal{U}} \setminus \mathcal{U}\right) \times \mathcal{U}$, where
$$\Phi(\zeta, z) = \langle P(\zeta, z), \zeta - z \rangle = \sum_{i=1}^n P_i(\zeta, z) \cdot (\zeta_i - z_i).$$

According to (1) we may assume that $\mathcal{U} = \widetilde{\mathcal{U}} \cap \mathbf{M}$ is a set of common zeros of smooth functions $\{\rho_k, \ k = 1, \ldots, m\}$. Then, using the q-concavity of $\mathbf{M}$ and applying Kohn's lemma to the set of functions $\{\rho_k\}$ we can construct a new set of functions $\tilde{\rho}_1, \ldots, \tilde{\rho}_m$ of the form:
$$\tilde{\rho}_k(z) = \rho_k(z) + A \cdot \left(\sum_{i=1}^m \rho_i^2(z)\right),$$
with large enough constant $A > 0$ and such that for any $z \in \mathbf{M}$ there exist an open neighborhood $\mathcal{U} \ni z$ and a family $E_{q+m}(\theta, z)$ of $q+m$ dimensional complex linear subspaces in $\mathbb{C}^n$ smoothly depending on $(\theta, z) \in \mathbb{S}^{m-1} \times \mathcal{U}$ and such that $-\mathcal{L}_z \tilde{\rho}_\theta$ is strictly negative on $E_{q+m}(\theta, z)$ with all negative eigenvalues not exceeding some $c < 0$.

To simplify notations we will assume that the functions $\rho_1, \ldots, \rho_m$ already satisfy this condition.

Let $E_{n-q-m}^\perp(\theta, z)$ be the family of $n-q-m$ dimensional subspaces in $T(\mathbf{G})$ orthogonal to $E_{q+m}(\theta, z)$ and let
$$a_j(\theta, z) = (a_{j1}(\theta, z), \ldots, a_{jn}(\theta, z)) \text{ for } j = 1, \ldots, n-q-m$$
be a set of $C^2$ smooth vector functions representing an orthonormal basis in $E_{n-q-m}^\perp(\theta, z)$.

Defining for $(\theta, z, w) \in \mathbb{S}^{m-1} \times \mathcal{U} \times \mathbb{C}^n$
$$A_j(\theta, z, w) = \sum_{i=1}^n a_{ji}(\theta, z) \cdot w_i, \quad (j = 1, \ldots, n-q-m)$$
we construct the form
$$\mathcal{A}(\theta, z, w) = \sum_{j=1}^{n-q-m} A_j(\theta, z, w) \cdot \bar{A}_j(\theta, z, w)$$

such that the hermitian form
$$\mathcal{L}_z \rho_\theta(w) + \mathcal{A}(\theta, z, w)$$
is strictly positive definite in $w$ for $(\theta, z) \in \mathbb{S}^{m-1} \times \mathcal{U}$.

Then we define for $\zeta, z \in \left(\widetilde{\mathcal{U}} \setminus \mathcal{U}\right) \times \mathcal{U}$:

$$\begin{aligned}
&Q_i^{(k)}(\zeta, z) = -\frac{\partial \rho_k}{\partial \zeta_i}(z) - \frac{1}{2} \sum_{j=1}^n \frac{\partial^2 \rho_k}{\partial \zeta_i \partial \zeta_j}(z)(\zeta_j - z_j), \\
&F^{(k)}(\zeta, z) = \langle Q^{(k)}(\zeta, z), \zeta - z \rangle, \\
&P_i(\zeta, z) = \sum_{k=1}^m \theta_k(\zeta) \cdot Q_i^{(k)}(\zeta, z) + \sum_{j=1}^{n-q-m} a_{ji}(\theta(\zeta), z) \cdot \bar{A}_j(\theta(\zeta), z, \zeta - z), \\
&\Phi(\zeta, z) = \langle P(\zeta, z), \zeta - z \rangle = \sum_{k=1}^m \theta_k(\zeta) \cdot F^{(k)}(\zeta, z) + \mathcal{A}(\theta(\zeta), z, \zeta - z)
\end{aligned} \tag{10}$$

with
$$\theta_k(\zeta) = -\frac{\rho_k(\zeta)}{\rho(\zeta)} \quad \text{for } k = 1, \ldots, m.$$

To prove that $P_i(\zeta, z)$ is a strong $\mathbf{M}$-barrier for some $U \ni z$ we consider the Taylor expansion of $\rho_k$ for $k = 1, \ldots, m$:
$$\rho_k(\zeta) = \rho_k(z) - 2\mathrm{Re}F^{(k)}(\zeta, z) + \mathcal{L}_z \rho_k(\zeta - z) + O(|\zeta - z|^3).$$

Then we obtain for some $U$ and $(\zeta, z) \in (U \setminus (U \cap \mathbf{M})) \times (U \cap \mathbf{M})$:
$$\mathrm{Re}\Phi(\zeta, z) = \sum_{k=1}^m \theta_k(\zeta) \cdot \mathrm{Re}F^{(k)}(\zeta, z) + \sum_{j=1}^{n-q-m} A_j(\theta(\zeta), \zeta, z) \cdot \bar{A}_j(\theta(\zeta), \zeta, z) \tag{11}$$
$$= \rho(\zeta) + \mathcal{L}_z \rho_\theta(\zeta - z) + \mathcal{A}(\theta(\zeta), z, \zeta - z) + O(|\zeta - z|^3),$$
which implies the existence of an open neighborhood $U \ni z$ in $\mathbb{C}^n$, satisfying (9).

Denoting $A_j(\zeta, z) := A_j(\theta(\zeta), z, \zeta - z)$ we obtain the following equalities that will be used in the further estimates
$$\bar{\partial}_\zeta \bar{A}_j(\zeta, z) = \mu_\tau^{(j)}(\zeta, z) + \mu_\nu^{(j)}(\zeta, z) \tag{12}$$
where
$$\mu_\tau^{(j)}(\zeta, z) = \sum_{i=1}^n \bar{a}_{ji}(\theta(\zeta), z) d\bar{\zeta}_i, \quad \mu_\nu^{(j)}(\zeta, z) = \sum_{i=1}^n (\bar{\zeta}_i - \bar{z}_i) \bar{\partial}_\zeta \bar{a}_{ji}(\theta(\zeta), z).$$

In our description of local integral formulas on $\mathbf{M}$ and in the future estimates we will also need the following notations.

We define the tubular neighborhood $\mathbf{G}_\epsilon$ of $\mathbf{M}$ in $\mathbf{G}$ as follows:
$$\mathbf{G}_\epsilon = \{z \in \mathbf{G} : \rho(z) < \epsilon\},$$
where $\rho(z) = \left(\sum_{k=1}^m \rho_k^2(z)\right)^{\frac{1}{2}}$. The boundary of $\mathbf{G}_\epsilon$ - $\mathbf{M}_\epsilon$ is defined by the condition
$$\mathbf{M}_\epsilon = \{z \in \mathbf{G} : \rho(z) = \epsilon\}.$$

For a sufficiently small neighborhood $\widetilde{\mathcal{U}} \in \mathbf{G}$ we may assume that functions
$$\rho_k(\zeta), \ \mathrm{Im}F^{(k)}(\zeta, z) \ \{k = 1, \ldots, m\}$$
have a nonzero jacobian with respect to $\mathrm{Re}\zeta_{i_1}, \ldots, \mathrm{Re}\zeta_{i_m}, \mathrm{Im}\zeta_{i_1}, \ldots, \mathrm{Im}\zeta_{i_m}$ for $z, \zeta \in \widetilde{\mathcal{U}}$. Therefore, for any fixed $z \in \widetilde{\mathcal{U}}$ these functions may be chosen as local $C^\infty$ coordinates in $\zeta$. We may also complement the functions above by holomorphic functions $w_j(\zeta) = u_j(\zeta) + iv_j(\zeta)$ with $j = 1, \ldots, n - m$ so that the functions
$$\rho_k(\zeta), \ \mathrm{Im}F^{(k)}(\zeta, z) \ \{k = 1, \ldots, m\},$$

$$u_j(\zeta), v_j(\zeta) \ \{j = 1, \ldots, n - m\},$$

represent a complete system of local coordinates in $\zeta \in \widetilde{\mathcal{U}}$ for any fixed $z \in \widetilde{\mathcal{U}}$.

The following complex valued vector fields on $\mathcal{U}$ for any fixed $z \in \mathcal{U}$

$$Y_{i,\zeta}(z) = \frac{\partial}{\partial \mathrm{Im} F^{(i)}(\zeta, z)} \text{ for } i = 1, \ldots, m,$$

$$W_{i,\zeta} = \frac{\partial}{\partial w_i}, \ \overline{W}_{i,\zeta} = \frac{\partial}{\partial \bar{w}_i} \text{ for } i = 1, \ldots, n - m,$$

represent a basis in $\mathbf{C}T(\mathbf{M})$ with vector fields $W_{i,\zeta}(z)$ being a basis in $T'(\mathbf{M})$ and $\overline{W}_{i,\zeta}(z)$ a basis in $T''(\mathbf{M})$. We denote

$$Y_\zeta(z) = -\sum_{k=1}^m \frac{\rho_k(\zeta)}{\rho(\zeta)} Y_{k,\zeta}(z).$$

We need also to consider local extensions of functions and forms from $\mathcal{U} = \widetilde{\mathcal{U}} \cap \mathbf{M}$ to $\widetilde{\mathcal{U}}$. To define appropriate functional spaces on these neighborhoods we consider fibration of $\widetilde{\mathcal{U}}$ by the manifolds $\widetilde{\mathcal{U}} \cap \mathbf{M}(\delta_1, \ldots, \delta_m)$ where

$$\mathbf{M}(\delta_1, \ldots, \delta_m) = \{z \in \mathbf{G} : \rho_1(z) = \delta_1, \ldots, \rho_m(z) = \delta_m\}.$$

Then we define spaces $\Pi^a\left(\{\rho\}, \widetilde{\mathcal{U}}(\epsilon)\right)$ of functions on $\widetilde{\mathcal{U}}(\epsilon) = \mathbf{G}_\epsilon \cap \widetilde{\mathcal{U}}$ in the same way as $\Pi^a(\mathcal{U})$ using the distribution $T^c(\mathbf{M}(\delta_1, \ldots, \delta_m))$ in $T(\mathbf{G})$. Namely, we define

$$\|h\|_{\Pi^a(\mathbf{M})} = \sup_{2k+s \leq p} \|\overbrace{D^c}^s \circ \overbrace{D}^k h\|_{\Gamma^\alpha(\mathbf{M})} + \sup_{2k+s \leq p-1} \|\overbrace{D^c}^s \circ \overbrace{D}^k h\|_{\Gamma^{1+\alpha}(\mathbf{M})}$$

with $D_i^c \in \mathbf{C}T^c(\mathbf{M}(\delta_1, \ldots, \delta_m))$ and $\|D_i\|_{C^{p+1}(\mathbf{M})}, \|D_i^c\|_{C^{p+1}(\mathbf{M})} \leq 1$ in the first term and $\|D_i\|_{C^p(\mathbf{M})}, \|D_i^c\|_{C^p(\mathbf{M})} \leq 1$ in the second term.

For a differential form $f = \sum_{I,J} f_{I,J}(z) dz^I \wedge d\bar{z}^J$ with $|I| = l$ and $|J| = r$ we say that $f \in \Pi_{(l,r)}^a\left(\{\rho\}, \widetilde{\mathcal{U}}(\epsilon)\right)$ if $f_{I,J} \in \Pi^a\left(\{\rho\}, \widetilde{\mathcal{U}}(\epsilon)\right)$.

We introduce a local extension operator for $\widetilde{\mathcal{U}}(\epsilon)$ and $\mathcal{U} = \widetilde{\mathcal{U}}(\epsilon) \cap \mathbf{M}$

$$E_\mathcal{U} : \Pi_{(0,r)}^a(\mathcal{U}) \to \Pi_{(0,r)}^a\left(\{\rho\}, \widetilde{\mathcal{U}}(\epsilon)\right)$$

which we define by extending all the coefficients of the differential form identically with respect to $\rho_1, \ldots, \rho_m$ in $\widetilde{\mathcal{U}}$. From the construction it follows that $E_\mathcal{U}$ satisfies the following estimate

$$\|E_\mathcal{U}(g)\|_{\Pi^a\left(\{\rho\}, \widetilde{\mathcal{U}}(\epsilon)\right)} \leq C \cdot \|g\|_{\Pi^a(\mathcal{U})}.$$

The following proposition provides local integral formula for $\bar{\partial}_\mathbf{M}$.

**Proposition 2.1.** *Let $\mathbf{M} \subset \mathbf{G}$ be a $C^\infty$ regular q-concave CR submanifold of the form (1) and let $\widetilde{\mathcal{U}}$ be an open neighborhood in $\mathbf{G}$ with analytic coordinates $z_1, \ldots, z_n$.*

*Then for $r = 1, \ldots, q - 1$ and any differential form $g \in C_{(0,r)}^\infty(\mathbf{M})$ with compact support in $\mathcal{U}$ the following equality*

$$g = \bar{\partial}_\mathbf{M} R_r(g) + R_{r+1}(\bar{\partial}_\mathbf{M} g) + H_r(g), \tag{13}$$

*holds, where*

$$R_r(g)(z)$$
$$= (-1)^r \frac{(n-1)!}{(2\pi i)^n} \cdot pr_\mathbf{M} \circ \lim_{\epsilon \to 0} \int_{\mathbf{M}_\epsilon \times [0,1]} \widetilde{g}(\zeta) \wedge \omega'_{r-1}\left((1-t)\frac{\bar{\zeta} - \bar{z}}{|\zeta - z|^2} + t\frac{P(\zeta, z)}{\Phi(\zeta, z)}\right) \wedge \omega(\zeta),$$

$$H_r(g)(z) = (-1)^r \frac{(n-1)!}{(2\pi i)^n} \cdot pr_\mathbf{M} \circ \lim_{\epsilon \to 0} \int_{\mathbf{M}_\epsilon} \widetilde{g}(\zeta) \wedge \omega'_r\left(\frac{P(\zeta, z)}{\Phi(\zeta, z)}\right) \wedge \omega(\zeta),$$

$\widetilde{g} = E_\mathcal{U}(g)$ is the extension of $g$, $\Phi(\zeta, z)$ is a local barrier for $\widetilde{\mathcal{U}}$ constructed in (10) and $pr_\mathbf{M}$ denotes the operator of projection to the space of tangential differential forms on $\mathbf{M}$.

We omit the proof of proposition 2.1 because it is completely analogous to the proof of formula (13) for another barrier function in [P].

Given above definitions of spaces $\Pi^a\left(\{\rho\}, \widetilde{\mathcal{U}}(\epsilon)\right)$ and of the extension operator $E_\mathcal{U}$ depend on the choice of functions $\rho_1, \ldots, \rho_m$. But we notice (cf. [P]) that the operators $R_r$ and $H_r$ are independent of the choice of functions $\rho_1, \ldots, \rho_m$ and of extension operator $E_\mathcal{U}$.

To construct now global formula on $\mathbf{M}$ we consider two finite coverings $\{\widetilde{\mathcal{U}}_\iota \subset \widetilde{\mathcal{U}}'_\iota\}$ of $\mathbf{G}$ and two partitions of unity $\{\vartheta_\iota\}$ and $\{\vartheta'_\iota\}$ subordinate to these coverings and such that $\vartheta'_\iota(z) = 1$ for $z \in \mathrm{supp}(\vartheta_\iota)$.

Applying proposition 2.1 to the form $\vartheta_\iota g$ in $\mathcal{U}'_\iota$ we obtain

$$\vartheta_\iota(z)g(z) = \bar{\partial}_\mathbf{M} R^\iota_r(\vartheta_\iota g)(z) + R^\iota_{r+1}(\bar{\partial}_\mathbf{M} \vartheta_\iota g)(z) + H^\iota_r(\vartheta_\iota g)(z).$$

Multiplying the equality above by $\vartheta'_\iota(z)$ and using equalities

$$\vartheta'_\iota(z) \cdot \bar{\partial}_\mathbf{M} R^\iota_r(\vartheta_\iota g)(z) = \bar{\partial}_\mathbf{M}\left[\vartheta'_\iota(z) \cdot R^\iota_r(\vartheta_\iota g)(z)\right] - \bar{\partial}_\mathbf{M}\vartheta'_\iota(z) \wedge R^\iota_r(\vartheta_\iota g)(z)$$

and

$$R^\iota_{r+1}(\bar{\partial}_\mathbf{M}\vartheta_\iota g)(z) = R^\iota_{r+1}(\bar{\partial}_\mathbf{M}\vartheta_\iota \wedge g)(z) + R^\iota_{r+1}(\vartheta_\iota \bar{\partial}_\mathbf{M} g)(z)$$

we obtain

$$\vartheta_\iota(z)g(z) = \bar{\partial}_\mathbf{M}\mathbf{R}^\iota_r(g)(z) + \mathbf{R}^\iota_{r+1}(\bar{\partial}_\mathbf{M} g)(z) + \mathbf{H}^\iota_r(g)(z) \tag{14}$$

with

$$\mathbf{R}^\iota_r(g)(z) = \vartheta'_\iota(z) \cdot R^\iota_r(\vartheta_\iota g)(z)$$

and

$$\mathbf{H}^\iota_r(g)(z) = -\bar{\partial}_\mathbf{M}\vartheta'_\iota(z) \wedge R^\iota_r(\vartheta_\iota g)(z) + \vartheta'_\iota(z) \cdot R^\iota_{r+1}(\bar{\partial}_\mathbf{M}\vartheta_\iota \wedge g)(z) + \vartheta'_\iota(z) \cdot H^\iota_r(\vartheta_\iota g)(z).$$

Adding equalities (14) for all $\iota$ we obtain

**Proposition 2.2.** *Let $\mathbf{M} \subset \mathbf{G}$ be a $C^\infty$ regular q-concave compact CR submanifold of the form (1).*

*Then for $r = 1, ..., q-1$ and any differential form $g \in C^\infty_{(0,r)}(\mathbf{M})$ the following equality*

$$g = \bar{\partial}_\mathbf{M}\mathbf{R}_r(g) + \mathbf{R}_{r+1}(\bar{\partial}_\mathbf{M} g) + \mathbf{H}_r(g), \tag{15}$$

*holds, where*

$$\mathbf{R}_r(g)(z) = \sum_\iota \vartheta'_\iota(z) \cdot R^\iota_r(\vartheta_\iota g)(z)$$

*and*

$$\mathbf{H}_r(g)(z) = \sum_\iota \left[-\bar{\partial}_\mathbf{M}\vartheta'_\iota(z) \wedge R^\iota_r(\vartheta_\iota g)(z) + R^\iota_{r+1}(\bar{\partial}_\mathbf{M}\vartheta_\iota \wedge g)(z) + \vartheta'_\iota(z) \cdot H^\iota_r(\vartheta_\iota g)(z)\right].$$

## 3. Boundedness of $\mathbf{R}_r$.

From the construction of operator $\mathbf{R}_r$ we conclude that in order to prove necessary estimates for operator $\mathbf{R}_r$ it suffices to prove these estimates for operator $R_r$. In the proposition below we state necessary estimates for operator $R_r$.

**Proposition 3.1.** *Let $a = p + \alpha$ with $0 < \alpha < 1$, $\mathbf{M} \subset \mathbf{G}$ be a $C^\infty$ regular $q$-concave CR submanifold of the form (1) and let $g \in \Pi^a_{(0,r)}(\mathbf{M})$ be a form with compact support in $\mathcal{U} = \widetilde{\mathcal{U}} \cap \mathbf{M}$. Then operator $R_r$, defined in (13) satisfies the following estimate*

$$\| R_r(g) \|_{\Pi^{a+1}_{(0,r-1)}(\mathcal{U})} < C \cdot \| g \|_{\Pi^a_{(0,r)}(\mathcal{U})} \tag{16}$$

*with a constant $C$ independent of $g$.*

In our proof of proposition 3.1 we will use the approximation of $R_r$ by the operators

$$R_r(\epsilon)(f)(z) = (-1)^r \cdot \mathrm{pr}_\mathbf{M} \circ \frac{(n-1)!}{(2\pi i)^n} \tag{17}$$

$$\times \sum_\iota \int_{\mathbf{M}_\epsilon \times [0,1]} \vartheta_\iota(\zeta) \widetilde{f}(\zeta) \wedge \omega'_{r-1}\left((1-t)\frac{\bar\zeta - \bar z}{|\zeta - z|^2} + t\frac{P^\iota(\zeta, z)}{\Phi^\iota(\zeta, z)}\right) \wedge \omega(\zeta)$$

when $\epsilon$ goes to 0.

Using equality (12) we obtain the following representation of kernels of these integrals on $\widetilde{\mathcal{U}} \times [0,1] \times \mathbf{M}$:

$$\vartheta_\iota(\zeta) \cdot \omega'_{r-1}\left((1-t)\frac{\bar\zeta - \bar z}{|\zeta - z|^2} + t\frac{P^\iota(\zeta, z)}{\Phi^\iota(\zeta, z)}\right) \wedge \omega(\zeta)\bigg|_{\widetilde{\mathcal{U}} \times [0,1] \times \mathbf{M}} \tag{18}$$

$$= \sum_{i,J} a_{(i,J)}(t, \zeta, z) dt \wedge \lambda^{i,J}_{r-1}(\zeta, z) + \sum_{i,J} b_{(i,J)}(t, \zeta, z) dt \wedge \gamma^{i,J}_{r-1}(\zeta, z),$$

where $i$ is an index, $J = \cup_{i=1}^8 J_i$ is a multiindex such that $i \notin J$, $a_{(i,J)}(t, \zeta, z)$ and $b_{(i,J)}(t, \zeta, z)$ are polynomials in $t$ with coefficients that are smooth functions of $z, \zeta$ and $\theta(\zeta)$, and $\lambda^{i,J}_{r-1}(\zeta, z)$ and $\gamma^{i,J}_{r-1}(\zeta, z)$ are defined as follows:

$$\lambda^{i,J}_{r-1}(\zeta, z) = \frac{1}{|\zeta - z|^{2(|J_1|+|J_5|+1)} \cdot \Phi(\zeta, z)^{n-|J_1|-|J_5|-1}} \tag{19}$$

$$\times \sum \mathrm{Det}\left[\bar\zeta - \bar z,\ Q^{(i)},\ \overbrace{d\bar\zeta}^{j \in J_1},\ \overbrace{\bar A \cdot \bar\partial_\zeta a}^{j \in J_2},\ \overbrace{a \cdot \mu_\nu}^{j \in J_3},\ \overbrace{a \cdot \mu_\tau}^{j \in J_4},\ \overbrace{d\bar z}^{j \in J_5},\ \overbrace{\bar A \cdot \bar\partial_z a}^{j \in J_6},\ \overbrace{a \cdot \bar\partial_z \bar A}^{j \in J_7},\ \overbrace{\bar\partial_z Q}^{j \in J_8}\right] \wedge \omega(\zeta),$$

and

$$\gamma^{i,J}_{r-1}(\zeta, z) = \frac{1}{|\zeta - z|^{2(|J_1|+|J_5|+1)} \cdot \Phi(\zeta, z)^{n-|J_1|-|J_5|-1}} \tag{20}$$

$$\times \sum \mathrm{Det}\left[\bar\zeta - \bar z,\ a_i \bar A_i,\ \overbrace{d\bar\zeta}^{j \in J_1},\ \overbrace{\bar A \cdot \bar\partial_\zeta a}^{j \in J_2},\ \overbrace{a \cdot \mu_\nu}^{j \in J_3},\ \overbrace{a \cdot \mu_\tau}^{j \in J_4},\ \overbrace{d\bar z}^{j \in J_5},\ \overbrace{\bar A \cdot \bar\partial_z a}^{j \in J_6},\ \overbrace{a \cdot \bar\partial_z \bar A}^{j \in J_7},\ \overbrace{\bar\partial_z Q}^{j \in J_8}\right] \wedge \omega(\zeta).$$

In the proof of the boundedness of operators $R_r$ we will need to know the differentiability properties of integrals with kernels $\lambda^{i,J}_{r-1}(\zeta, z)$ and $\gamma^{i,J}_{r-1}(\zeta, z)$. In the lemmas below we prepare necessary tools.

We introduce kernels

$$\mathcal{K}^I_{d,h}(\zeta,z) = \frac{\{\rho(\zeta)\}^{I_1}(\zeta-z)^{I_2}(\bar\zeta-\bar z)^{I_3}}{|\zeta-z|^d \cdot \Phi(\zeta,z)^h} \wedge \overbrace{d\rho_i}^{i\in I_4} \wedge \overbrace{d\theta_i(\zeta)}^{i\in I_5} \wedge d\sigma_{2n-m}(\zeta),$$

where $I = \cup_{j=1}^5 I_j$ and $I_j$ for $j=1,\dots,5$ are multiindices such that $I_1$ contains $m$ indices, $I_2$, $I_3$ contain $n$ indices, $I_4 \cup I_5$ contains $m-1$ indices, $|I_4|+|I_5| = m-1$, and $\{\rho(\zeta)\}^{I_1} = \prod_{i_s \in I_1} \rho_s(\zeta)^{i_s}$, $(\zeta-z)^{I_2} = \prod_{i_s \in I_2}(\zeta_s - z_s)^{i_s}$, $(\bar\zeta - \bar z)^{I_3} = \prod_{i_s \in I_3}(\bar\zeta_s - \bar z_s)^{i_s}$.

For kernels $\mathcal{K}^I_{d,h}$ we introduce the following notation

$$k\left(\mathcal{K}^I_{d,h}\right) = d - |I_2| - |I_3|,$$
$$h\left(\mathcal{K}^I_{d,h}\right) = h,$$
$$l\left(\mathcal{K}^I_{d,h}\right) = |I_1| + |I_4|.$$

In two lemmas below we describe smoothing properties of kernels $\mathcal{K}^I_{d,h}$.

**Lemma 3.2.** *Let $\mathcal{U} = \widetilde{\mathcal{U}} \cap \mathbf{M}$ be a neighborhood with $\Phi(\zeta,z)$ constructed in (10), $\mathcal{U}(\epsilon) = \widetilde{\mathcal{U}} \cap \mathbf{M}_\epsilon$, and let $g(\zeta,z,\theta,t)$ be a smooth form with compact support in $\widetilde{\mathcal{U}}_\zeta \times \widetilde{\mathcal{U}}_z \times \mathbb{S}^{n-1} \times [0,1]$.*

*Then for $g(\zeta,z,t) = g(\zeta,z,\theta(\zeta),t)$ and a vector field*

$$D_z = \sum_{j=1}^n a_j(z)\frac{\partial}{\partial z_j} + \sum_{j=1}^n b_j(z)\frac{\partial}{\partial \bar z_j} \in \mathbf{C}T(\mathbf{G})$$

*such that $D|_\mathcal{U} \in \mathbf{C}T(\mathbf{M})$ the following equality holds*

$$D_z\left(\int_{\mathcal{U}(\epsilon)\times[0,1]} g(\zeta,z,t) \cdot \mathcal{K}^I_{d,h}(\zeta,z)dt\right) \tag{21}$$

$$= \int_{\mathcal{U}(\epsilon)\times[0,1]} [D_z g(\zeta,z,t)] \cdot \mathcal{K}^I_{d,h}(\zeta,z)dt + \int_{\mathcal{U}(\epsilon)\times[0,1]} [D_\zeta g(\zeta,z,t)] \cdot \mathcal{K}^I_{d,h}(\zeta,z)dt$$

$$+ \sum_{S,a,b} \int_{\mathcal{U}(\epsilon)\times[0,1]} c_{\{S,a,b\}}(\zeta,z,t) \cdot g(\zeta,z,t) \cdot \mathcal{K}^S_{a,b}(\zeta,z)dt$$

$$+ \sum_{L,i} \int_{\mathcal{U}(\epsilon)\times[0,1]} c_{\{L,i\}}(\zeta,z,t) \cdot [Y_\zeta(z)g(\zeta,z,t)] \cdot \mathcal{K}^L_{i,h}(\zeta,z)dt,$$

*where $c_{\{S,a,b\}}(\zeta,z,t), c_{\{L,i,e,k\}}(\zeta,z,t)$ are $C^\infty$ functions of $\zeta, z, \theta(\zeta), t$, vector field $D_\zeta$ is defined as*

$$D_\zeta = \sum_{j=1}^n a_j(\zeta)\frac{\partial}{\partial \zeta_j} + \sum_{j=1}^n b_j(\zeta)\frac{\partial}{\partial \bar\zeta_j}$$

*and kernels $\mathcal{K}^S_{a,b}$ and $\mathcal{K}^L_{i,h}$ satisfy the following conditions*

$$k\left(\mathcal{K}^S_{a,b}\right) + b - l\left(\mathcal{K}^S_{a,b}\right) \le k\left(\mathcal{K}^I_{d,h}\right) + h - l\left(\mathcal{K}^I_{d,h}\right),$$
$$k\left(\mathcal{K}^S_{a,b}\right) + 2b - 2l\left(\mathcal{K}^S_{a,b}\right) \le k\left(\mathcal{K}^I_{d,h}\right) + 2h - 2l\left(\mathcal{K}^I_{d,h}\right), \tag{22}$$

*and*

$$k\left(\mathcal{K}^L_{i,h}\right) - l\left(\mathcal{K}^L_{i,h}\right) + 1 \le k\left(\mathcal{K}^I_{d,h}\right) - l\left(\mathcal{K}^I_{d,h}\right),$$
$$k\left(\mathcal{K}^L_{i,h}\right) - 2l\left(\mathcal{K}^L_{i,h}\right) + 1 \le k\left(\mathcal{K}^I_{d,h}\right) - 2l\left(\mathcal{K}^I_{d,h}\right). \tag{23}$$

**Proof.** To prove the lemma we represent the integral from the left hand side of (21) as

$$D_z \left( \int_{\mathcal{U}(\epsilon) \times [0,1]} g(\zeta, z, t) \cdot \mathcal{K}^I_{d,h}(\zeta, z) dt \right) \tag{24}$$

$$= \int_{\mathcal{U}(\epsilon) \times [0,1]} [D_z g(\zeta, z, t)] \cdot \mathcal{K}^I_{d,h}(\zeta, z) dt - \int_{\mathcal{U}(\epsilon) \times [0,1]} g(\zeta, z, t) \left[ D_\zeta \mathcal{K}^I_{d,h}(\zeta, z) \right] dt$$

$$+ \int_{\mathcal{U}(\epsilon) \times [0,1]} g(\zeta, z, t) \left[ (D_z + D_\zeta) \mathcal{K}^I_{d,h}(\zeta, z) \right] dt.$$

To transform the second term of the right hand side of (24) we apply integration by parts and obtain

$$\int_{\mathcal{U}(\epsilon) \times [0,1]} g(\zeta, z, t) \left[ D_\zeta \mathcal{K}^I_{d,h}(\zeta, z) \right] dt = - \int_{\mathcal{U}(\epsilon) \times [0,1]} [D_\zeta g(\zeta, z, t)] \mathcal{K}^I_{d,h}(\zeta, z) dt$$

$$+ \sum_{S,a,b} \int_{\mathcal{U}(\epsilon) \times [0,1]} c_{\{S,a,b\}}(\zeta, z, t) \cdot g(\zeta, z, t) \cdot \mathcal{K}^S_{a,b}(\zeta, z) dt$$

with kernels $\mathcal{K}^S_{a,b}$ satisfying (22).

To transform the third term of the right hand side of (24) we will use the estimates below that follow from the definitions of $F^{(k)}(\zeta, z)$ and $\mathcal{A}(\zeta, z)$ and from the fact that $D_z|_{\mathcal{U}} \in \mathbf{CT}(\mathbf{M})$

$$(D_z + D_\zeta) \mathcal{A}(\zeta, z) = \mathcal{O}\left(|\zeta - z|^2\right),$$

$$(D_z + D_\zeta) \operatorname{Re} F^{(k)}(\zeta, z) = \mathcal{O}\left(|\zeta - z|^2\right),$$

$$(D_z + D_\zeta) \operatorname{Im} F^{(k)}(\zeta, z) = \mathcal{O}\left(|\zeta - z|\right), \tag{25}$$

$$(D_z + D_\zeta)(\zeta_j - z_j) = \mathcal{O}\left(|\zeta - z|\right),$$

$$(D_z + D_\zeta)(\bar{\zeta}_j - \bar{z}_j) = \mathcal{O}\left(|\zeta - z|\right).$$

Applying operators $D_z$ and $D_\zeta$ to $\mathcal{K}^I_{d,h}(\zeta, z)$ and using estimates (25) we obtain

$$(D_z + D_\zeta) \mathcal{K}^I_{d,h}(\zeta, z) \tag{26}$$

$$= (-h) \frac{\{\rho(\zeta)\}^{I_1}(\zeta - z)^{I_2}(\bar{\zeta} - \bar{z})^{I_3}}{|\zeta - z|^d \cdot \Phi(\zeta, z)^{h+1}} \left( [(D_z + D_\zeta) \mathcal{A}] \overbrace{\wedge d\rho_i}^{i \in I_4} \wedge \overbrace{d\theta_i(\zeta)}^{i \in I_5} \wedge d\sigma_{2n-m}(\zeta) \right.$$

$$+ \sum_{k=1}^m \theta_k(\zeta) \cdot \left[ (D_z + D_\zeta) \operatorname{Re} F^{(k)} \right] \overbrace{\wedge d\rho_i}^{i \in I_4} \wedge \overbrace{d\theta_i(\zeta)}^{i \in I_5} \wedge d\sigma_{2n-m}(\zeta)$$

$$\left. + \sum_{k=1}^m \theta_k(\zeta) \cdot \left[ (D_z + D_\zeta) \operatorname{Im} F^{(k)} \right] \overbrace{\wedge d\rho_i}^{i \in I_4} \wedge \overbrace{d\theta_i(\zeta)}^{i \in I_5} \wedge d\sigma_{2n-m}(\zeta) \right)$$

$$+ \sum_{S,a,b} c_{\{S,a,b\}}(\zeta, z, t) \cdot \mathcal{K}^S_{a,b}(\zeta, z)$$

with kernels $\mathcal{K}^S_{a,b}$ satisfying (22).

From estimates (25) we conclude that the first two terms of the right hand side of (26) can be represented as linear combinations with $C^\infty \left( \widetilde{\mathcal{U}}_\zeta \times \widetilde{\mathcal{U}}_z \times \mathbb{S}^{n-1} \times [0,1] \right)$ coefficients of kernels $\mathcal{K}^S_{d,h+1}$ with

$$|S_2| + |S_3| = |I_2| + |I_3| + 2,$$

and therefore
$$k\left(\mathcal{K}^S_{d,h+1}\right) \le k\left(\mathcal{K}^I_{d,h+1}\right) - 2.$$
Then corresponding integrals can be represented in the form
$$\sum_S \int_{\mathcal{U}(\epsilon)\times[0,1]} c_{\{S\}}(\zeta,z,t) \cdot g(\zeta,z,t) \cdot \mathcal{K}^S_{d,h+1}(\zeta,z) dt$$
with kernels $\mathcal{K}^S_{d,h+1}$ satisfying (22).

The third term of the right hand side of (26) we represent as
$$\sum_{k=1}^m \frac{\theta_k(\zeta) \cdot \left[(D_z+D_\zeta)\operatorname{Im} F^{(k)}(\zeta,z)\right]\{\rho(\zeta)\}^{I_1}(\zeta-z)^{I_2}(\bar\zeta-\bar z)^{I_3}}{|\zeta-z|^d}$$
$$\times d_\zeta\left[\frac{\overbrace{\wedge d\rho_i}^{i\in I_4}\wedge \overbrace{d\theta_i(\zeta)}^{i\in I_5}\wedge (d_\zeta\operatorname{Im}\Phi(\zeta,z)\lrcorner\, d\sigma_{2n-m}(\zeta))}{\Phi(\zeta,z)^h}\right]$$
$$+ \sum_{S,a,b} c_{\{S\}}(\zeta,z,t)\cdot \mathcal{K}^S_{a,b}(\zeta,z)$$
with kernels $\mathcal{K}^S_{a,b}$ satisfying (22).

Applying integration by parts to the corresponding integrals and using the third estimate from (25) we obtain
$$\sum_{k=1}^m \int_{\mathcal{U}(\epsilon)\times[0,1]} g(\zeta,z,t)\cdot \frac{\theta_k(\zeta)\cdot \left[(D_z+D_\zeta)\operatorname{Im} F^{(k)}(\zeta,z)\right]\{\rho(\zeta)\}^{I_1}(\zeta-z)^{I_2}(\bar\zeta-\bar z)^{I_3}}{|\zeta-z|^d}$$
$$\times d_\zeta\left[\frac{\overbrace{\wedge d\rho_i}^{i\in I_4}\wedge \overbrace{d\theta_i(\zeta)}^{i\in I_5}\wedge (d_\zeta\operatorname{Im}\Phi(\zeta,z)\lrcorner\, d\sigma_{2n-m}(\zeta))}{\Phi(\zeta,z)^h}\right]$$
$$= \sum_{k=1}^m \int_{\mathcal{U}(\epsilon)\times[0,1]} \theta_k(\zeta) g(\zeta,z,t)\cdot Y_\zeta(z)\left[\frac{\left[(D_z+D_\zeta)\operatorname{Im} F^{(k)}(\zeta,z)\right]\{\rho(\zeta)\}^{I_1}(\zeta-z)^{I_2}(\bar\zeta-\bar z)^{I_3}}{|\zeta-z|^d}\right]$$
$$\times \frac{\overbrace{\wedge d\rho_i}^{i\in I_4}\wedge \overbrace{d\theta_i(\zeta)}^{i\in I_5}\wedge (d_\zeta\operatorname{Im}\Phi(\zeta,z)\lrcorner\, d\sigma_{2n-m}(\zeta))}{\Phi(\zeta,z)^h}$$
$$+ \sum_{k=1}^m \int_{\mathcal{U}(\epsilon)\times[0,1]} \theta_k(\zeta)\left[Y_\zeta(z) g(\zeta,z,t)\right]\cdot \frac{\left[(D_z+D_\zeta)\operatorname{Im} F^{(k)}(\zeta,z)\right]\{\rho(\zeta)\}^{I_1}(\zeta-z)^{I_2}(\bar\zeta-\bar z)^{I_3}}{|\zeta-z|^d}$$
$$\times \frac{\overbrace{\wedge d\rho_i}^{i\in I_4}\wedge \overbrace{d\theta_i(\zeta)}^{i\in I_5}\wedge (d_\zeta\operatorname{Im}\Phi(\zeta,z)\lrcorner\, d\sigma_{2n-m}(\zeta))}{\Phi(\zeta,z)^h}$$
$$= \sum_{S,a} \int_{\mathcal{U}(\epsilon)\times[0,1]} c_{\{S,a\}}(\zeta,z,t)\cdot g(\zeta,z,t)\cdot \mathcal{K}^S_{a,h}(\zeta,z) dt$$
$$+ \sum_{L,i} \int_{\mathcal{U}(\epsilon)\times[0,1]} c_{\{L,i\}}(\zeta,z,t)\cdot \left[Y_\zeta(z) g(\zeta,z,t)\right]\cdot \mathcal{K}^L_{i,h}(\zeta,z) dt$$
with kernels satisfying (22) and (23). $\square$

**Lemma 3.3.** *Let* $\mathcal{U} = \widetilde{\mathcal{U}} \cap \mathbf{M}$ *be a neighborhood with* $\Phi(\zeta, z)$ *constructed in (10),* $\mathcal{U}(\epsilon) = \widetilde{\mathcal{U}} \cap \mathbf{M}_\epsilon$ *and let* $g(\zeta, z, \theta, t)$ *be a smooth form with compact support in* $\widetilde{\mathcal{U}}_\zeta \times \widetilde{\mathcal{U}}_z \times \mathbb{S}^{n-1} \times [0, 1]$.

*Then for* $g(\zeta, z, t) = g(\zeta, z, \theta(\zeta), t)$ *and vector fields*

$$D^c_{i,z} = \sum_{j=1}^n a_{ij}(z)\frac{\partial}{\partial z_j} + \sum_{j=1}^n b_{ij}(z)\frac{\partial}{\partial \bar{z}_j} \in \mathbf{CT}(\mathbf{G}) \quad (i = 1, \ldots, k)$$

*such that* $D|_{\mathcal{U}} \in \mathbf{CT}^c(\mathbf{M})$ *the following equality holds*

$$D^c_{k,z} \circ \cdots \circ D^c_{1,z} \left( \int_{\mathcal{U}(\epsilon) \times [0,1]} g(\zeta, z, t) \cdot \mathcal{K}^I_{d,h}(\zeta, z) dt \right) \tag{27}$$

$$= \sum_{\|R\| \leq k} \sum_{\{R,S,a,b\}} \int_{\mathcal{U}(\epsilon) \times [0,1]} c_{\{R,S,a,b\}}(\zeta, z, t) \cdot \left[ \{Y_\zeta(z), D^c_\zeta, D^c_z\}^R g(\zeta, z, t) \right] \cdot \mathcal{K}^S_{a,b}(\zeta, z) dt$$

$$+ \sum_{\|R\|=k+1} \sum_{\{R,L,e,m\}} \int_{\mathcal{U}(\epsilon) \times [0,1]} c_{\{R,L,e,m\}}(\zeta, z, t) \cdot \left[ \{Y_\zeta(z), D^c_\zeta, D^c_z\}^R g(\zeta, z, t) \right] \cdot \mathcal{K}^L_{e,m}(\zeta, z) dt,$$

*where* $R = (r_1, r_2, r_3)$, $\{Y_\zeta(z), D^c_\zeta, D^c_z\}^R$ *denotes a composition of* $r_1$ *differentiations* $Y_\zeta(z)$, $r_2$ *differentiations* $D^c_{i,\zeta}$ *and* $r_3$ *differentiations* $D^c_{i,z}$ *with* $\|R\| = 2r_1 + r_2 + r_3$,

$$c_{\{R,S,a,b\}}(\zeta, z, \theta, t), c_{\{R,L,e,m\}}(\zeta, z, \theta, t) \in C^\infty\left(\widetilde{\mathcal{U}}_\zeta \times \widetilde{\mathcal{U}}_z \times \mathbb{S}^{n-1} \times [0,1]\right),$$

*kernels* $\mathcal{K}^S_{a,b}$ *satisfy inequalities (22) and kernels* $\mathcal{K}^L_{e,m}$ *satisfy inequalities (23).*

**Proof.** We prove the lemma by induction with respect to $k$. The statement of the lemma for $k=1$ is a straightforward corollary of lemma 3.2.

To prove the step of induction we assume that equality (27) holds for $k$. To prove it for $k+1$ we apply $D^c_{k+1,z}$ to the right hand side of (27). To apply $D^c_{k+1,z}$ to a term with a kernel $\mathcal{K}^L_{e,m}$ we simply apply the differentiation to the kernel. Using estimate

$$\left| D^c_{k+1,z} \mathrm{Im} F^{(i)}(\zeta, z) \right| = \mathcal{O}\left(|\zeta - z|\right) \quad \text{for} \quad i = 1, \ldots, m,$$

we conclude that after differentiation $D^c_{k+1,z}$ of the kernel $\mathcal{K}^L_{e,m}$ the only possible changes in indices are

(i) $k\left(\mathcal{K}^L_{e,m}\right)$ increases by 1,

(ii) $h\left(\mathcal{K}^L_{e,m}\right)$ increases by 1 and $k\left(\mathcal{K}^L_{e,m}\right)$ decreases by 1.

Therefore, since kernel $\mathcal{K}^L_{e,m}$ satisfies conditions (23) we obtain

$$D^c_{k+1,z} \mathcal{K}^L_{e,m}(\zeta, z)$$

$$= \sum_{\|R\| \leq k+1} \sum_{\{R,S,a,b\}} c_{\{R,S,a,b\}}(\zeta, z, t) \cdot \left[ \{Y_\zeta(z), D^c_\zeta, D^c_z\}^R g(\zeta, z, t) \right] \cdot \mathcal{K}^S_{a,b}(\zeta, z)$$

with kernels $\mathcal{K}^S_{a,b}$ satisfying conditions (22).

To apply $D^c_{k+1,z}$ to a term with kernel $\mathcal{K}^S_{a,b}$ we use lemma 3.2 and obtain

$$D^c_{k+1,z} \int_{\mathcal{U}(\epsilon) \times [0,1]} c_{\{R,S,a,b\}}(\zeta, z, t) \cdot \left[ \{Y_\zeta(z), D^c_\zeta, D^c_z\}^R g(\zeta, z, t) \right] \cdot \mathcal{K}^S_{a,b}(\zeta, z) dt \tag{28}$$

$$= \sum_{\|P\| \leq k+1} \sum_{\{P,I,d,j\}} \int_{\mathcal{U}(\epsilon) \times [0,1]} c_{\{P,I,d,j\}}(\zeta, z, t) \cdot \left[ \{Y_\zeta(z), D^c_\zeta, D^c_z\}^P g(\zeta, z, t) \right] \cdot \mathcal{K}^I_{d,h}(\zeta, z) dt$$

$$+ \sum_{\|P\|=k+2} \sum_{\{P,L,e,m\}} \int_{\mathcal{U}(\epsilon) \times [0,1]} c_{\{P,L,e,m\}}(\zeta, z, t) \cdot \left[ \{Y_\zeta(z), D^c_\zeta, D^c_z\}^P g(\zeta, z, t) \right] \cdot \mathcal{K}^L_{e,m}(\zeta, z) dt,$$

with kernels $\mathcal{K}^I_{d,h}$ satisfying (22) and kernels $\mathcal{K}^L_{e,m}$ satisfying (23). $\square$

The following two simple lemmas will be used in the further estimates.

**Lemma 3.4.** *Let $\mathbf{M}$ be a generic CR submanifold in the unit ball $\mathbf{B}^n$ in $\mathbb{C}^n$ of the form:*
$$\mathbf{M} = \{z \in \mathbf{B}^n : \rho_1(z) = \cdots = \rho_m(z) = 0\},$$
*where $\{\rho_k\}$, $k=1,\ldots,m$ ($m < n$) are real valued functions of the class $C^\infty$ satisfying*
$$\partial \rho_1 \wedge \cdots \wedge \partial \rho_m \neq 0 \quad \text{on } \mathbf{M}.$$
*Then for any point $\zeta_0 \in \mathbf{M}$ there exists a neighborhood $\mathbf{V}_\epsilon(\zeta_0) = \{\zeta : |\zeta - \zeta_0| < \epsilon\}$ such that for any $n \geq s > n-m$ and $p > 2n-s-m$ the following representation holds in $\mathbf{V}_\epsilon$:*
$$d\bar\zeta_{i_1} \wedge \ldots d\bar\zeta_{i_p} \wedge d\zeta_{k_1} \wedge \ldots d\zeta_{k_s} = \sum d\rho_{j_1} \wedge \ldots d\rho_{j_{p-(2n-s-m)}} \wedge g^{i_1\ldots i_p}_{j_1\ldots j_{p-(2n-s-m)}}(\zeta) \tag{29}$$
*with $g^{i_1\ldots i_p}_{j_1\ldots j_{p-(2n-s-m)}}$ of the class $C^\infty(\mathbf{V}_\epsilon)$.* $\square$

**Lemma 3.5.** *Let*
$$\mathbf{B}(\delta) = \{(\eta, w) \in \mathbb{R}^m \times \mathbb{C}^{n-m} : \sum_{i=1}^m \eta_i^2 + \sum_{i=1}^{n-m} |w|^2 < \delta\},$$
$$\mathbf{V}(\delta) = \{(\eta, w) \in \mathbb{R}^m \times \mathbb{C}^{n-m} : |\eta_1| + \sum_{i=2}^m |\eta_i|^2 + \sum_{i=1}^{n-m} |w|^2 < \delta^2\},$$
$$K\{\alpha, k, h\}(\eta, w, \epsilon)$$
$$= \frac{\wedge_{i=1}^m d\eta_i \wedge_{i=1}^{n-m}(dw_i \wedge d\bar w_i)}{(\epsilon + \sum_{i=1}^m |\eta_i| + \sum_{i=1}^{n-m} |w_i|)^k (\sqrt{\epsilon} + \sqrt{|\eta_1|} + \sum_{i=2}^m |\eta_i| + \sum_{i=1}^{n-m} |w_i|)^{2h-\alpha}},$$
*with $0 \leq \alpha < 1$ and $k, 2h \in \mathbb{Z}$.*

*Let*
$$\mathcal{I}_1\{\alpha, k, h\}(\epsilon, \delta) = \int_{\mathbf{V}(\delta)} K\{\alpha, k, h\}(\eta, w, \epsilon),$$
*and*
$$\mathcal{I}_2\{\alpha, k, h\}(\epsilon, \delta) = \int_{\mathbf{B}(1) \setminus \mathbf{V}(\delta)} K\{\alpha, k, h\}(\eta, w, \epsilon).$$
*Then*
$$\mathcal{I}_1\{0, k, h\}(\epsilon, \delta)$$
$$= \begin{cases} \mathcal{O}\left(\epsilon^{2n-m-k-h} \cdot (\log \epsilon)^2\right) & \text{if } k \geq 2n-m-1 \text{ and } k+h \geq 2n-m, \\ \mathcal{O}(\delta) & \text{if } k \geq 2n-m-1 \text{ and } k+h \leq 2n-m-1, \\ \mathcal{O}\left(\epsilon^{(2n-m-k-2h+1)/2} \cdot \log \epsilon\right) & \text{if } k \leq 2n-m-2 \text{ and } k+2h \geq 2n-m+1, \\ \mathcal{O}(\delta) & \text{if } k \leq 2n-m-2 \text{ and } k+2h \leq 2n-m, \end{cases}$$

$$\mathcal{I}_1\{\alpha, k, h\}(\epsilon, \delta)$$
$$= \mathcal{O}(\delta^\alpha) \text{ if } \alpha > 0, \text{ and } \begin{cases} \text{if } k \geq 2n-m-1 \text{ and } k+h \leq 2n-m-1/2, \\ \text{if } k \leq 2n-m-2 \text{ and } k+2h \leq 2n-m+1, \end{cases}$$

$$\mathcal{I}_2\{0, k, h\}(\epsilon, \delta)$$

$$= \begin{cases} \mathcal{O}(1) & \begin{cases} \text{if } k \geq 2n-m-1 \text{ and } k+h \leq 2n-m-1, \\ \text{if } k \leq 2n-m-2 \text{ and } k+2h \leq 2n-m, \end{cases} \\ \mathcal{O}(\log \delta) & \text{if } k \leq 2n-m-2 \text{ and } k+2h \leq 2n-m+1, \end{cases}$$

$$\mathcal{I}_2\{\alpha, k, h\}(\epsilon, \delta)$$
$$= \mathcal{O}\left(\delta^{\alpha-1}\right) \text{ if } \alpha \geq 0, \text{ and } \begin{cases} \text{if } k \geq 2n-m-1 \text{ and } k+h \leq 2n-m, \\ \text{if } k \leq 2n-m-2 \text{ and } k+2h \leq 2n-m+2, \end{cases}$$

$$\mathcal{I}_2\{\alpha, k, h\}(\epsilon, \delta)$$
$$= \mathcal{O}\left(\delta^{\alpha-2}\right) \text{ if } \alpha \geq 0, \text{ and } \begin{cases} \text{if } k \geq 2n-m-1 \text{ and } k+h \leq 2n-m+1/2, \\ \text{if } k \leq 2n-m-2 \text{ and } k+2h \leq 2n-m+3. \end{cases}$$

**Proof.** Introducing $W = \left(\sum_{i=1}^{n-m} |w_i|^2\right)^{\frac{1}{2}}$, $\eta = \left(\sum_{i=2}^m \eta_i^2\right)^{\frac{1}{2}}$ and $V = (W^2 + \eta^2)^{\frac{1}{2}}$ we obtain:

$$K\{\alpha, k, h\}(\eta, w, \epsilon) = \frac{\eta^{m-2} W^{2n-2m-1} d\eta_1 \wedge d\eta \wedge dW}{(\epsilon + \eta_1 + \eta + W)^k \left(\sqrt{\epsilon} + \sqrt{\eta_1} + \eta + W\right)^{2h-\alpha}}$$
$$= \frac{V^{2n-m-2} d\eta_1 \wedge dV}{(\epsilon + \eta_1 + V)^k \left(\sqrt{\epsilon} + \sqrt{\eta_1} + V\right)^{2h-\alpha}}.$$

□

**Proof of proposition 3.1.**

According to (17) in order to prove the statement of the proposition it suffices to prove the estimates

$$\|\int_{\mathcal{U}(\epsilon) \times [0,1]} a_{(i,J)}(t, \zeta, z) dt \wedge \widetilde{g}(\zeta) \wedge \lambda_{r-1}^{i,J}(\zeta, z)\|_{\Pi^{a+1}(\mathcal{U})} \leq C \cdot \|g\|_{\Pi^a(\mathcal{U})}, \tag{30}$$

$$\|\int_{\mathcal{U}(\epsilon) \times [0,1]} b_{(i,J)}(t, \zeta, z) dt \wedge \widetilde{g}(\zeta) \wedge \gamma_{r-1}^{i,J}(\zeta, z)\|_{\Pi^{a+1}(\mathcal{U})} \leq C \cdot \|g\|_{\Pi^a(\mathcal{U})}$$

with constant $C$ independent of $g$ and $\epsilon$.

Using estimates
$$|\bar{A}| = \mathcal{O}(|\zeta - z|), \quad |\mu_\nu| = \mathcal{O}(|\zeta - z|) \tag{31}$$
for the terms of determinants in (19) and (20) and applying lemma 3.4 to the differential form

$$\overbrace{d\bar{\zeta}}^{|J_1|+r} \wedge \overbrace{\mu_\tau}^{|J_4|} \wedge \omega(\zeta)$$

we obtain representations

$$a_{(i,J)}(t, \zeta, z) dt \wedge \widetilde{g}(\zeta) \wedge \lambda_{r-1}^{i,J}(\zeta, z) = \sum_{\{I,d,j\}} c_{\{I,d,j\}}(\zeta, z, t) \widetilde{g}(\zeta) \mathcal{K}_{d,h}^I(\zeta, z), \tag{32}$$

and

$$b_{(i,J)}(t, \zeta, z) dt \wedge \widetilde{g}(\zeta) \wedge \gamma_{r-1}^{i,J}(\zeta, z) = \sum_{\{I,d,j\}} c_{\{I,d,j\}}(\zeta, z, t) \widetilde{g}(\zeta) \mathcal{K}_{d,h}^I(\zeta, z). \tag{33}$$

Multiindices $I_i$ for $i = 1, \ldots, 5$ and indices $d, h$ in (32) satisfy the conditions

$$\begin{aligned} d &= 2(|J_1| + |J_5| + 1), \\ h &= n - |J_1| - |J_5| - 1, \\ |I_2| + |I_3| &= 1 + |J_2| + |J_3| + |J_6|, \\ |I_4| &= |J_1| + |J_4| + r + m - n. \end{aligned} \tag{34}$$

Multiindices $I_i$ for $i = 1, \ldots, 5$ and indices $d, h$ in (33) satisfy the conditions

$$\begin{aligned} d &= 2(|J_1| + |J_5| + 1), \\ h &= n - |J_1| - |J_5| - 1, \\ |I_2| + |I_3| &= 2 + |J_2| + |J_3| + |J_6|, \\ |I_4| &= |J_1| + |J_4| + r + m - n. \end{aligned} \qquad (35)$$

Using representations (32) and (33) we reduce the problem of proving (30) to each term

$$\widetilde{g}(\zeta) \cdot c_{\{I,d,j\}}(\zeta, z, t) \mathcal{K}^I_{d,h}(\zeta, z)$$

of the right hand sides of these representations.

To further reduce the proof of (30) using lemmas 3.2 and 3.3 we will prove that only the case $0 < a < 2$ has to be considered.

**Lemma 3.6.** *Statement of proposition 3.1 follows from the corresponding statement for $0 < a < 2$. Namely, for the proof of estimates (30) it suffices to prove that for $0 < \alpha < 1$*

$$\begin{aligned} \left\| \int_{\mathcal{U}(\epsilon) \times [0,1]} c(\zeta, z, t) \widetilde{g}(\zeta) \mathcal{K}^I_{d,h}(\zeta, z) dt \right\|_{\Pi^{1+\alpha}(\mathcal{U})} &\leq C \cdot \|g\|_{\Pi^\alpha(\mathcal{U})}, \\ \left\| \int_{\mathcal{U}(\epsilon) \times [0,1]} c(\zeta, z, t) \widetilde{g}(\zeta) \mathcal{K}^I_{d,h}(\zeta, z) dt \right\|_{\Pi^{2+\alpha}(\mathcal{U})} &\leq C \cdot \|g\|_{\Pi^{1+\alpha}(\mathcal{U})} \end{aligned} \qquad (36)$$

*for the kernels $\mathcal{K}^I_{d,h}$ obtained from $\lambda^{i,J}_{r-1}$ and $\gamma^{i,J}_{r-1}$ after application of lemmas 3.2 and 3.3.*

**Proof.** Let us fix $a = p + \alpha$. As follows from the definition of spaces $\Pi^{a+1}$ in order to prove proposition 3.1 we have to prove that

(a) for any set of tangent vector fields $D_1, \ldots, D_k, D^c_1, \ldots, D^c_s$ on $\mathbf{M}$ such that

$$\|D_i\|_{C^{p+2}(\mathbf{M})}, \|D^c_i\|_{C^{p+2}(\mathbf{M})} \leq 1$$

with $D^c_i \in \mathbf{C}T^c(\mathbf{M})$ and $2k + s \leq p + 1$

$$\| \overbrace{D^c}^{s} \circ \overbrace{D}^{k} \int_{\mathcal{U}(\epsilon) \times [0,1]} \widetilde{g}(\zeta) \cdot c_{\{I,d,j\}}(\zeta, z, t) \mathcal{K}^I_{d,h}(\zeta, z) dt \|_{\Gamma^\alpha(\mathbf{M})} < C \|g\|_{\Pi^a(\mathbf{M})},$$

(b) for any set of tangent vector fields $D_1, \ldots, D_k, D^c_1, \ldots, D^c_s$ on $\mathbf{M}$ such that

$$\|D_i\|_{C^{p+1}(\mathbf{M})}, \|D^c_i\|_{C^{p+1}(\mathbf{M})} \leq 1$$

with $D^c_i \in \mathbf{C}T^c(\mathbf{M})$ and $2k + s \leq p$

$$\| \overbrace{D^c}^{s} \circ \overbrace{D}^{k} \int_{\mathcal{U}(\epsilon) \times [0,1]} \widetilde{g}(\zeta) \cdot c_{\{I,d,j\}}(\zeta, z, t) \mathcal{K}^I_{d,h}(\zeta, z) dt \|_{\Gamma^{1+\alpha}(\mathbf{M})} < C \|g\|_{\Pi^a(\mathbf{M})}$$

with constant $C$ independent of $g$.

Let us fix $k, s \in \mathbb{Z}^+ \cup \{0\}$ such that $2k + s = p + 1$ and $k, s \neq 0$. At first we apply an operator

$$\overbrace{D^c}^{s-1} \circ \overbrace{D}^{k-1}$$

to an integral

$$\int_{\mathcal{U}(\epsilon) \times [0,1]} \widetilde{g}(\zeta) \cdot c_{\{I,d,j\}}(\zeta, z, t) \mathcal{K}^I_{d,h}(\zeta, z) dt$$

using lemmas 3.2 and 3.3.

Then we obtain a representation

$$\overbrace{D^c}^{s-1} \circ \overbrace{D}^{k-1} \int_{\mathcal{U}(\epsilon) \times [0,1]} \widetilde{g}(\zeta) \cdot c_{\{I,d,j\}}(\zeta, z, t) \mathcal{K}^I_{d,h}(\zeta, z) dt \quad (37)$$

$$= \sum_{\|P\| \leq p-2} \sum_{\{P,S,a,b\}} \int_{\mathcal{U}(\epsilon) \times [0,1]} c_{\{P,S,a,b\}}(\zeta, z, t) \cdot \left[\{Y_\zeta(z), D_\zeta, D^c_\zeta\}^P \widetilde{g}(\zeta)\right] \cdot \mathcal{K}^S_{a,b}(\zeta, z) dt$$

$$+ \sum_{\|R\|=p-1} \sum_{\{R,L,e,m\}} \int_{\mathcal{U}(\epsilon) \times [0,1]} c_{\{R,L,e,m\}}(\zeta, z, t) \cdot \left[\{Y_\zeta(z), D_\zeta, D^c_\zeta\}^R \widetilde{g}(\zeta)\right] \cdot \mathcal{K}^L_{e,m}(\zeta, z) dt,$$

with the same notation as in (27) and kernels $\mathcal{K}^S_{a,b}$ and $\mathcal{K}^L_{e,m}$ satisfying (22) and (23) respectively.

We notice that according to the assumption of proposition 3.1 we have

$$\begin{aligned} \{Y_\zeta(z), D_\zeta, D^c_\zeta\}^P \widetilde{g} &\in \Pi^{2+\alpha}\left(\{\rho\}, \widetilde{\mathcal{U}}(\epsilon)\right), \\ \{Y_\zeta(z), D_\zeta, D^c_\zeta\}^R \widetilde{g} &\in \Pi^{1+\alpha}\left(\{\rho\}, \widetilde{\mathcal{U}}(\epsilon)\right). \end{aligned} \quad (38)$$

We apply a differentiation $D_z$ once more to terms in the first sum of (37) using lemma 3.2 and represent the result as

$$D_z \left[ \sum_{\|P\| \leq p-2} \sum_{\{P,S,a,b\}} \int_{\mathcal{U}(\epsilon) \times [0,1]} c_{\{P,S,a,b\}}(\zeta, z, t) \cdot \left[\{Y_\zeta(z), D_\zeta, D^c_\zeta\}^P \widetilde{g}(\zeta)\right] \cdot \mathcal{K}^S_{a,b}(\zeta, z) dt \right]$$

$$= \sum_{\|P\| \leq p} \sum_{\{P,I,d,j\}} \int_{\mathcal{U}(\epsilon) \times [0,1]} c_{\{P,I,d,j\}}(\zeta, z, t) \cdot \left[\{Y_\zeta(z), D_\zeta, D^c_\zeta\}^P \widetilde{g}(\zeta)\right] \cdot \mathcal{K}^I_{d,h}(\zeta, z) dt \quad (39)$$

with kernels $\mathcal{K}^I_{d,h}$ satisfying (22) and $\{Y_\zeta(z), D_\zeta, D^c_\zeta\}^P \widetilde{g} \in \Pi^\alpha\left(\{\rho\}, \widetilde{\mathcal{U}}(\epsilon)\right)$.

We apply a differentiation $D^c_z$ once more to terms of the second sum of (37) differentiating kernels $\mathcal{K}^L_{e,m}$ and using the same arguments as in the proof of lemma 3.3 we represent the result as

$$D^c_z \left[ \sum_{\|R\|=p-1} \sum_{\{R,L,e,m\}} \int_{\mathcal{U}(\epsilon) \times [0,1]} c_{\{R,L,e,m\}}(\zeta, z, t) \cdot \left[\{Y_\zeta(z), D_\zeta, D^c_\zeta\}^R \widetilde{g}(\zeta)\right] \cdot \mathcal{K}^L_{e,m}(\zeta, z) dt \right]$$

$$= \sum_{\|P\| \leq p-1} \sum_{\{P,I,d,j\}} \int_{\mathcal{U}(\epsilon) \times [0,1]} c_{\{P,I,d,j\}}(\zeta, z, t) \cdot \left[\{Y_\zeta(z), D_\zeta, D^c_\zeta\}^P \widetilde{g}(\zeta)\right] \cdot \mathcal{K}^I_{d,h}(\zeta, z) dt \quad (40)$$

with kernels $\mathcal{K}^I_{d,h}$ satisfying (22) and $\{Y_\zeta(z), D_\zeta, D^c_\zeta\}^P \widetilde{g} \in \Pi^{1+\alpha}\left(\{\rho\}, \widetilde{\mathcal{U}}(\epsilon)\right)$.

From (37), (38), (39) and (40) we conclude that in order to prove the statement of proposition 3.1 for $k, s \in \mathbb{Z}^+$ such that $2k + s = p + 1$ it suffices to prove (36) for the kernels $\mathcal{K}^I_{d,h}$ obtained from $\lambda^{i,J}_{r-1}$ and $\gamma^{i,J}_{r-1}$ after application of lemmas 3.2 and 3.3.

Analogous arguments show that in the case $s = 0$ or $k = 0$ and $2k + s = p + 1$ estimates (36) are also sufficient. Namely, in the case $s = 0$ after application of lemma 3.2 $k - 1$ times we obtain the right hand side of equality (40) and in the case $k = 0$ we use lemma 3.3 for $s - 2$ differentiations $D^c_z$ and then equality (40) for kernels $\mathcal{K}^L_{e,m}$.

For $k, s \in \mathbb{Z}^+ \cup \{0\}$ such that $2k + s = p$ and $s \neq 0$ we apply an operator

$$\overbrace{D^c}^{s-1} \circ \overbrace{D}^{k}$$

to an integral
$$\int_{\mathcal{U}(\epsilon)\times[0,1]} \widetilde{g}(\zeta) \cdot c_{\{I,d,j\}}(\zeta,z,t)\mathcal{K}^I_{d,h}(\zeta,z)dt$$
and use lemmas 3.2 and 3.3.

Then we obtain a representation
$$\overbrace{D^c}^{s-1} \circ \overbrace{D}^{k} \int_{\mathcal{U}(\epsilon)\times[0,1]} \widetilde{g}(\zeta) \cdot c_{\{I,d,j\}}(\zeta,z,t)\mathcal{K}^I_{d,h}(\zeta,z)dt \tag{41}$$

$$= \sum_{\|P\|\leq p-1} \sum_{\{P,S,a,b\}} \int_{\mathcal{U}(\epsilon)\times[0,1]} c_{\{P,S,a,b\}}(\zeta,z,t) \cdot \left[\{Y_\zeta(z), D_\zeta, D^c_\zeta\}^P \widetilde{g}(\zeta)\right] \cdot \mathcal{K}^S_{a,b}(\zeta,z)dt$$

$$+ \sum_{\|R\|=p} \sum_{\{R,L,e,m\}} \int_{\mathcal{U}(\epsilon)\times[0,1]} c_{\{R,L,e,m\}}(\zeta,z,t) \cdot \left[\{Y_\zeta(z), D_\zeta, D^c_\zeta\}^R \widetilde{g}(\zeta)\right] \cdot \mathcal{K}^L_{e,m}(\zeta,z)dt,$$

with kernels $\mathcal{K}^S_{a,b}$ and $\mathcal{K}^L_{e,m}$ satisfying (22) and (23) respectively, and
$$\{Y_\zeta(z), D_\zeta, D^c_\zeta\}^P \widetilde{g} \in \Pi^{1+\alpha}\left(\{\rho\}, \widetilde{\mathcal{U}}(\epsilon)\right),$$
$$\{Y_\zeta(z), D_\zeta, D^c_\zeta\}^R \widetilde{g} \in \Pi^\alpha\left(\{\rho\}, \widetilde{\mathcal{U}}(\epsilon)\right).$$

Applying one more differentiation $D^c$ to terms of the second sum of (41) by differentiating kernels $\mathcal{K}^L_{e,m}$ we conclude that in the case of $k, s \in \mathbb{Z}^+ \cup \{0\}$ such that $2k + s = p$ estimates (36) are also sufficient. $\square$

We will prove estimates (36) as a corollary of two lemmas below.

**Lemma 3.7.** *Let $0 < \alpha < 1$, $g \in \Gamma^\alpha(\mathcal{U})$ be a function with compact support and $\widetilde{g} \in \Gamma^\alpha\left(\{\rho\}, \widetilde{\mathcal{U}}(\epsilon_0)\right)$ its extension. Let $\mathcal{K}^S_{a,b}$ satisfy conditions*
$$\begin{aligned} k\left(\mathcal{K}^S_{a,b}\right) + b - l\left(\mathcal{K}^S_{a,b}\right) &\leq 2n - m - 2, \\ k\left(\mathcal{K}^S_{a,b}\right) + 2b - 2l\left(\mathcal{K}^S_{a,b}\right) &\leq 2n - m. \end{aligned} \tag{42}$$

*Then*
$$f_\epsilon(z) := \left(\int_{\mathcal{U}(\epsilon)\times[0,1]} c(\zeta,z,t)\widetilde{g}(\zeta)\mathcal{K}^S_{a,b}(\zeta,z)dt\right) \in \Pi^{1+\alpha}(\mathcal{U})$$
*and*
$$\|f_\epsilon\|_{\Pi^{1+\alpha}(\mathcal{U})} \leq C \cdot \|g\|_{\Gamma^\alpha(\mathcal{U})}$$
*with $C$ independent of $g$ and $\epsilon$.*

**Proof.** From the definition of spaces $\Pi^a$ we conclude that statement of the lemma would follow from inclusions
$$f_\epsilon \in \Lambda^{\frac{1+\alpha}{2}}(\mathcal{U}) \quad \text{and} \quad D^c f_\epsilon \in \Gamma^\alpha(\mathcal{U}).$$

We start with inclusion $f_\epsilon \in \Lambda^{\frac{1+\alpha}{2}}(\mathcal{U})$. For $w \in \mathcal{U}$ and arbitrary $\delta > 0$ we introduce neighborhoods
$$W(w, \epsilon, \sqrt{\delta}) = \left\{\zeta \in \widetilde{\mathcal{U}} : \rho(\zeta) = \epsilon, \ |\Phi(\zeta, w)| \leq c \cdot \delta\right\},$$
such that for $|z - w| \leq \delta$
$$W(w, \epsilon, \sqrt{\delta}) \subset W(z, \epsilon, C\sqrt{\delta})$$

with constants $c, C > 0$ independent of $w$, $z$ and $\delta$.

Then we represent $f_\epsilon(z)$ for $z$ such that $|z - w| \leq \delta$ as

$$f_\epsilon(z) = g(w) \cdot \int_{\mathcal{U}(\epsilon) \times [0,1]} c(\zeta, z, t) \mathcal{K}_{a,b}^S(\zeta, z) dt \qquad (43)$$

$$+ \int_{W(w,\epsilon,\sqrt{\delta}) \times [0,1]} (\widetilde{g}(\zeta) - g(w)) \, c(\zeta, z, t) \mathcal{K}_{a,b}^S(\zeta, z) dt$$

$$+ \int_{(\mathcal{U}(\epsilon) \setminus W(w,\epsilon,\sqrt{\delta})) \times [0,1]} (\widetilde{g}(\zeta) - g(w)) \, c(\zeta, z, t) \mathcal{K}_{a,b}^S(\zeta, z) dt.$$

Applying lemma 3.2 to the first term of the right hand side of (43) we obtain

$$D_z \left( g(w) \int_{\mathcal{U}(\epsilon) \times [0,1]} c(\zeta, z, t) \cdot \mathcal{K}_{a,b}^S(\zeta, z) dt \right)$$

$$= g(w) \cdot \sum_{I,d,j} \int_{\mathcal{U}(\epsilon) \times [0,1]} c_{\{I,d,j\}}(\zeta, z, t) \cdot \mathcal{K}_{d,h}^I(\zeta, z) dt$$

with $c_{\{I,d,j\}}(\zeta, z, t) = c_{\{I,d,j\}}(\zeta, z, \theta(\zeta), t) \in C^\infty \left( \widetilde{\mathcal{U}}_\zeta \times \widetilde{\mathcal{U}}_z \times \mathbb{S}^{n-1} \times [0,1] \right)$, and kernels $\mathcal{K}_{d,h}^I$ satisfying (42).

Then applying formula

$$d\rho_i|_{\mathcal{U}(\epsilon)} = \epsilon d\theta_i$$

we obtain

$$\left| g(w) \cdot D_z \left( \int_{\mathcal{U}(\epsilon) \times [0,1]} c(\zeta, z, t) \mathcal{K}_{a,b}^S(\zeta, z) dt \right) \right|$$

$$= \|g\|_{\Lambda^0(\mathcal{U})} \cdot \mathcal{O} \left( \epsilon^{l(\mathcal{K})} \int_{\mathcal{U}(\epsilon) \times [0,1]} \frac{\overbrace{\wedge_i d\theta_i(\zeta)}^{m-1} \wedge d\sigma_{2n-m}(\zeta)}{|\zeta - z|^{k(\mathcal{K})} \cdot |\Phi(\zeta, z)|^{h(\mathcal{K})}} \right)$$

$$= \|g\|_{\Lambda^0(\mathcal{U})} \cdot \mathcal{O} \left( \epsilon^{l(\mathcal{K})} \cdot \mathcal{I}_1 \{0, k(\mathcal{K}), h(\mathcal{K})\} (\epsilon, 1) \right).$$

Applying lemma 3.5 to the last term of the estimate above we have

$$\epsilon^{l(\mathcal{K})} \cdot \mathcal{I}_1 \{0, k(\mathcal{K}), h(\mathcal{K})\} (\epsilon, 1) = \begin{cases} \mathcal{O} \left( \epsilon^{2n-m-k-h+l} \cdot (\log \epsilon)^2 \right) \\ \mathcal{O} \left( \epsilon^{(2n-m-k-2h+2l+1)/2} \cdot \log \epsilon \right) \\ \mathcal{O}(1) \end{cases} = \mathcal{O}(1),$$

which shows that the first term of the right hand side of (43) is in $\Lambda^1(\mathcal{U})$.

For the second term of the right hand side of (43) we have

$$\left| \int_{W(w,\epsilon,\sqrt{\delta}) \times [0,1]} (\widetilde{g}(\zeta) - g(w)) \, c(\zeta, z, t) \mathcal{K}_{a,b}^S(\zeta, z) dt \right|$$

$$= \|g\|_{\Gamma^\alpha(\mathcal{U})} \cdot \delta^{\frac{\alpha}{2}} \cdot \mathcal{O} \left( \epsilon^{l(\mathcal{K})} \cdot \int_{W(z,\epsilon,C\sqrt{\delta}) \times [0,1]} \frac{\overbrace{\wedge_i d\theta_i(\zeta)}^{m-1} \wedge d\sigma_{2n-m}(\zeta)}{|\zeta - z|^{k(\mathcal{K})} \cdot |\Phi(\zeta, z)|^{h(\mathcal{K})}} \right)$$

$$= \|g\|_{\Gamma^\alpha(\mathcal{U})} \cdot \delta^{\frac{\alpha}{2}} \cdot \mathcal{O} \left( \epsilon^{l(\mathcal{K})} \cdot \mathcal{I}_1 \{0, k(\mathcal{K}), h(\mathcal{K})\} \left( \epsilon, \sqrt{\delta} \right) \right) = \|g\|_{\Gamma^\alpha(\mathcal{U})} \cdot \mathcal{O} \left( \delta^{\frac{1+\alpha}{2}} \right),$$

where we again used lemma 3.5 and the estimate

$$|\widetilde{g}(\zeta) - g(w)| = \|g\|_{\Gamma^\alpha(\mathcal{U})} \cdot \mathcal{O} \left( \delta^{\frac{\alpha}{2}} \right)$$

for $\zeta \in W(z, \epsilon, C\sqrt{\delta})$.

For the third term of the right hand side of (43) using estimates

$$|\widetilde{g}(\zeta) - g(z)| = \|g\|_{\Gamma^\alpha(\mathcal{U})} \cdot \mathcal{O}\left(|\Phi(\zeta, z)|^{\alpha/2}\right)$$

and

$$\left|F^{(k)}(\zeta, z) - F^{(k)}(\zeta, w)\right| = \mathcal{O}(\delta) \tag{44}$$

and lemma 3.5 we obtain

$$\left|\int_{(\mathcal{U}(\epsilon) \setminus W(w,\epsilon,\sqrt{\delta})) \times [0,1]} (\widetilde{g}(\zeta) - g(w)) \, c(\zeta, z, t) \mathcal{K}^S_{a,b}(\zeta, z) dt\right.$$

$$\left. - \int_{(\mathcal{U}(\epsilon) \setminus W(w,\epsilon,\sqrt{\delta})) \times [0,1]} (\widetilde{g}(\zeta) - g(w)) \, c(\zeta, w, t) \mathcal{K}^S_{a,b}(\zeta, w) dt\right|$$

$$\leq \left|\int_{(\mathcal{U}(\epsilon) \setminus W(w,\epsilon,\sqrt{\delta})) \times [0,1]} (\widetilde{g}(\zeta) - g(w)) \, c(\zeta, z, t) \left[\mathcal{K}^S_{a,b}(\zeta, z) - \mathcal{K}^S_{a,b}(\zeta, w)\right] dt\right|$$

$$+ \left|\int_{(\mathcal{U}(\epsilon) \setminus W(w,\epsilon,\sqrt{\delta})) \times [0,1]} (\widetilde{g}(\zeta) - g(w)) \left[c(\zeta, z, t) - c(\zeta, w, t)\right] \mathcal{K}^S_{a,b}(\zeta, w) dt\right|$$

$$= \|g\|_{\Gamma^\alpha(\mathcal{U})} \cdot \delta \cdot \mathcal{O}\left[\mathcal{I}_2\left\{\alpha, k(\mathcal{K}) + 1, b - l(\mathcal{K})\right\}\left(\epsilon, \sqrt{\delta}\right)\right.$$

$$+ \mathcal{I}_2\left\{\alpha, k(\mathcal{K}), b - l(\mathcal{K}) + 1\right\}\left(\epsilon, \sqrt{\delta}\right)$$

$$\left. + \mathcal{I}_2\left\{\alpha, k(\mathcal{K}), b - l(\mathcal{K})\right\}\left(\epsilon, \sqrt{\delta}\right)\right] = \|g\|_{\Gamma^\alpha(\mathcal{U})} \cdot \mathcal{O}\left(\delta^{\frac{1+\alpha}{2}}\right).$$

Representation (43) together with the estimates above show that

$$\|f_\epsilon\|_{\Lambda^{\frac{1+\alpha}{2}}(\mathcal{U})} \leq C \cdot \|g\|_{\Gamma^\alpha(\mathcal{U})}$$

uniformly with respect to $\epsilon$.

To complete the proof of the lemma we have to prove that

$$\begin{aligned}\|D^c f_\epsilon\|_{\Lambda^\alpha_c(\mathcal{U})} &\leq C \cdot \|g\|_{\Gamma^\alpha(\mathcal{U})}, \\ \|D^c f_\epsilon\|_{\Lambda^{\frac{\alpha}{2}}(\mathcal{U})} &\leq C \cdot \|g\|_{\Gamma^\alpha(\mathcal{U})},\end{aligned} \tag{45}$$

where differentiation $D^c \in \mathbf{CT}^c(\mathbf{M})$ and

$$\|h\|_{\Lambda^\alpha_c(\mathcal{U})} = \sup\left\{|h(x(\cdot))|_{\Lambda^\alpha([0,1])}\right\}$$

with the *sup* taken over all curves $x : [0, 1] \to \mathbf{M}$ satisfying (2).

To prove these estimates we use the following representation

$$D^c_z f_\epsilon(z) = g(z) \cdot D^c_z \int_{\mathcal{U}(\epsilon) \times [0,1]} c(\zeta, z, t) \mathcal{K}^S_{a,b}(\zeta, z) dt \tag{46}$$

$$+ \int_{\mathcal{U}(\epsilon) \times [0,1]} (\widetilde{g}(\zeta) - g(z)) \, D^c_z \left[c(\zeta, z, t) \mathcal{K}^S_{a,b}(\zeta, z)\right] dt.$$

Applying lemma 3.3 to the first term of (46) and using lemma 3.5 we conclude that

$$D^c_z \int_{\mathcal{U}(\epsilon) \times [0,1]} c(\zeta, z, t) \mathcal{K}^S_{a,b}(\zeta, z) dt \in \Lambda^1(\mathcal{U})$$

and therefore the first term of the right hand side of (46) is in $\Gamma^\alpha(\mathcal{U})$.

Thus statement of the lemma is now reduced to the proof of estimates (45) for

$$\int_{\mathcal{U}(\epsilon) \times [0,1]} (\widetilde{g}(\zeta) - g(z)) \, D^c_z \left[c(\zeta, z, t) \mathcal{K}^S_{a,b}(\zeta, z)\right] dt. \tag{47}$$

We will prove the first estimate from (45) for (47). Proof of the second one is analogous with a few changes that we will point out.

For $w \in \mathcal{U}$ and arbitrary $\delta > 0$ we introduce neighborhoods

$$W(w, \epsilon, \delta) = \left\{ \zeta \in \widetilde{\mathcal{U}} : \rho(\zeta) = \epsilon, \ |\Phi(\zeta, w)| \leq c\delta^2 \right\},$$
$$W'(w, \epsilon, \delta) = \left\{ \zeta \in \widetilde{\mathcal{U}} : \rho(\zeta) = \epsilon, \ |\Phi(\zeta, w)| \leq c'\delta^2 \right\}, \tag{48}$$

with $c < c'$ and such that if $|z - w| \leq \delta$ and there exists a curve $x : [0, 1] \to \mathbf{M}$ satisfying (2) with $x(0) = w$ and $x(1) = z$ then

$$W'(w, \epsilon, \delta) \subset W(z, \epsilon, C\delta)$$

with $c, c', C > 0$ independent of $z, w, \delta$.

We also introduce a function $\phi_w(\zeta) \in C^\infty(\widetilde{\mathcal{U}} \setminus \mathcal{U})$ such that $0 \leq \phi_w(\zeta) \leq 1$ and

$$\phi_w \equiv 1 \ \text{on} \ W(w, \epsilon, \delta), \quad \phi_w \equiv 0 \ \text{on} \ \mathcal{U}(\epsilon) \setminus W'(w, \epsilon, \delta),$$
$$|D\phi_w(\zeta)| = \mathcal{O}(1/\delta^2), \quad |D^c\phi_w(\zeta)| = \mathcal{O}(1/\delta), \tag{49}$$

for vector fields $D, D^c \in \mathbf{CT}(\mathbf{G})$ such that $D|_\mathcal{U} \in \mathbf{CT}(\mathbf{M})$ and $D^c|_\mathcal{U} \in \mathbf{CT}^c(\mathbf{M})$.

Then for fixed $w \in \mathcal{U}$ and $z \in \mathcal{U}$ such that $|z - w| \leq \delta$ and such that there exists a curve $x : [0, 1] \to \mathbf{M}$ satisfying (2) with $x(0) = w$ and $x(1) = z$ we consider the following representation for integral in (47)

$$\int_{\mathcal{U}(\epsilon) \times [0,1]} (\widetilde{g}(\zeta) - g(z)) \, D_z^c \left[ c(\zeta, z, t) \mathcal{K}_{a,b}^S(\zeta, z) \right] dt \tag{50}$$

$$- \int_{\mathcal{U}(\epsilon) \times [0,1]} (\widetilde{g}(\zeta) - g(w)) \, D_w^c \left[ c(\zeta, w, t) \mathcal{K}_{a,b}^S(\zeta, w) \right] dt$$

$$= \int_{\mathcal{U}(\epsilon) \times [0,1]} (\widetilde{g}(\zeta) - g(z)) \, \phi_w(\zeta) D_z^c \left[ c(\zeta, z, t) \mathcal{K}_{a,b}^S(\zeta, z) \right] dt$$

$$- \int_{\mathcal{U}(\epsilon) \times [0,1]} (\widetilde{g}(\zeta) - g(w)) \, \phi_w(\zeta) D_w^c \left[ c(\zeta, w, t) \mathcal{K}_{a,b}^S(\zeta, w) \right] dt$$

$$+ \int_{\mathcal{U}(\epsilon) \times [0,1]} \left[ (\widetilde{g}(\zeta) - g(z)) (1 - \phi_w(\zeta)) \, D_z^c \left[ c(\zeta, z, t) \mathcal{K}_{a,b}^S(\zeta, z) \right] \right.$$

$$\left. - (\widetilde{g}(\zeta) - g(w)) (1 - \phi_w(\zeta)) \, D_w^c \left[ c(\zeta, w, t) \mathcal{K}_{a,b}^S(\zeta, w) \right] \right] dt.$$

To estimate integrals in the right hand side of (50) we use estimate

$$D_z^c \left[ F^{(k)}(\zeta, z) \right] = \mathcal{O}(|\zeta - z|) \tag{51}$$

for $k = 1, \ldots, m$ and obtain a representation

$$D_z^c \left[ c(\zeta, z, t) \mathcal{K}_{a,b}^S(\zeta, z) \right] = \sum_{I, d, j} c_{\{I, d, j\}}(\zeta, z, t) \cdot \mathcal{K}_{d,h}^I(\zeta, z) \tag{52}$$

with $c_{\{I,d,j\}}(\zeta, z, t) = c_{\{I,d,j\}}(\zeta, z, \theta(\zeta), t) \in C^\infty\left(\widetilde{\mathcal{U}}_\zeta \times \widetilde{\mathcal{U}}_z \times \mathbb{S}^{n-1} \times [0, 1]\right)$, and kernels $\mathcal{K}_{d,h}^I$ satisfying

$$k\left(\mathcal{K}_{d,h}^I\right) + h - l\left(\mathcal{K}_{d,h}^I\right) \leq 2n - m - 1,$$
$$k\left(\mathcal{K}_{d,h}^I\right) + 2h - 2l\left(\mathcal{K}_{d,h}^I\right) \leq 2n - m + 1. \tag{53}$$

Then for the first two integrals of the right hand side of (50) using representation (52) with kernels satisfying conditions (53) and lemma 3.5 we obtain

$$\left| \int_{\mathcal{U}(\epsilon) \times [0,1]} (\widetilde{g}(\zeta) - g(z)) \, \phi_w(\zeta) c(\zeta, z, t) \mathcal{K}_{d,h}^I(\zeta, z) dt \right|$$

$$= \mathcal{O}\left(\left|\int_{W'(w,\epsilon,\delta)\times[0,1]} (\tilde{g}(\zeta) - g(z))\, c(\zeta, z, t)\mathcal{K}^I_{d,h}(\zeta, z)\, dt\right|\right)$$

$$= \|g\|_{\Gamma^\alpha(\mathcal{U})} \cdot \mathcal{O}\left(\int_{W(z,\epsilon,C\delta)} \frac{\overbrace{\wedge_i d\theta_i(\zeta)}^{m-1} \wedge d\sigma_{2n-m}(\zeta)}{|\zeta - z|^{k(\mathcal{K})} \cdot |\Phi(\zeta, z)|^{h(\mathcal{K})-l(\mathcal{K})-\alpha/2}}\right)$$

$$= \|g\|_{\Gamma^\alpha(\mathcal{U})} \cdot \mathcal{O}\left(\mathcal{I}_1\{\alpha, k(\mathcal{K}), h(\mathcal{K}) - l(\mathcal{K})\}(\epsilon, \delta)\right) = \|g\|_{\Gamma^\alpha(\mathcal{U})} \cdot \mathcal{O}(\delta^\alpha).$$

To estimate the third integral of the right hand side of (50) we represent it as

$$\int_{\mathcal{U}(\epsilon)\times[0,1]} (\tilde{g}(\zeta) - g(z))(1 - \phi_w(\zeta))\, D^c_z\left[c(\zeta, z, t)\mathcal{K}^S_{a,b}(\zeta, z)\right] dt \tag{54}$$

$$-\int_{\mathcal{U}(\epsilon)\times[0,1]} (\tilde{g}(\zeta) - g(w))(1 - \phi_w(\zeta))\, D^c_w\left[c(\zeta, w, t)\mathcal{K}^S_{a,b}(\zeta, w)\right] dt$$

$$= \int_{\mathcal{U}(\epsilon)\times[0,1]} (\tilde{g}(\zeta) - g(z))(1 - \phi_w(\zeta))$$

$$\times \left(D^c_z\left[c(\zeta, z, t)\mathcal{K}^S_{a,b}(\zeta, z)\right] - D^c_w\left[c(\zeta, w, t)\mathcal{K}^S_{a,b}(\zeta, w)\right]\right) dt$$

$$+ [g(w) - g(z)] \cdot \int_{\mathcal{U}(\epsilon)\times[0,1]} (1 - \phi_w(\zeta))\, D^c_w\left[c(\zeta, w, t)\mathcal{K}^S_{a,b}(\zeta, w)\right] dt.$$

For the first integral of the right hand side of (54) we use representation (52) with kernels satisfying (53) and estimate

$$\left|F^{(k)}(\zeta, z) - F^{(k)}(\zeta, w)\right| = \mathcal{O}\left(\delta \cdot |\zeta - z| + \delta^2\right) \tag{55}$$

for $k = 1, \ldots, m$ and $\zeta \in \mathcal{U}(\epsilon) \setminus W(w, \epsilon, \delta)$.

Then using lemma 3.5 we obtain

$$\left|\int_{\mathcal{U}(\epsilon)\times[0,1]} (\tilde{g}(\zeta) - g(z))(1 - \phi_w(\zeta))\right.$$

$$\left.\times \left(D^c\left[c(\zeta, z, t)\mathcal{K}^S_{a,b}(\zeta, z)\right] - D^c\left[c(\zeta, w, t)\mathcal{K}^S_{a,b}(\zeta, w)\right]\right) dt\right|$$

$$= \left|\sum_{I,d,j} \int_{\mathcal{U}(\epsilon)\times[0,1]} (\tilde{g}(\zeta) - g(z))(1 - \phi_w(\zeta))\right.$$

$$\left.\times \left[c_{\{I,d,j\}}(\zeta, z, t) \cdot \mathcal{K}^I_{d,h}(\zeta, z) - c_{\{I,d,j\}}(\zeta, w, t) \cdot \mathcal{K}^I_{d,h}(\zeta, w)\right] dt\right|$$

$$= \|g\|_{\Gamma^\alpha(\mathcal{U})} \cdot \left[\delta \cdot \mathcal{O}\left(\mathcal{I}_2\{\alpha, k(\mathcal{K}), h(\mathcal{K}) - l(\mathcal{K})\}(\epsilon, \delta)\right)\right.$$

$$+ \delta \cdot \mathcal{O}\left(\mathcal{I}_2\{\alpha, k(\mathcal{K}) + 1, h(\mathcal{K}) - l(\mathcal{K})\}(\epsilon, \delta)\right)$$

$$+ \delta \cdot \mathcal{O}\left(\mathcal{I}_2\{\alpha, k(\mathcal{K}) - 1, h(\mathcal{K}) - l(\mathcal{K}) + 1\}(\epsilon, \delta)\right)$$

$$\left.+ \delta^2 \cdot \mathcal{O}\left(\mathcal{I}_2\{\alpha, k(\mathcal{K}), h(\mathcal{K}) - l(\mathcal{K}) + 1\}(\epsilon, \delta)\right)\right] = \|g\|_{\Gamma^\alpha(\mathcal{U})} \cdot \mathcal{O}(\delta^\alpha).$$

To obtain necessary estimate for the second integral of the right hand side of (54) it suffices to prove the following estimate

$$\left|\int_{\mathcal{U}(\epsilon)\times[0,1]} (1 - \phi_w(\zeta))\, D^c_w\left[c(\zeta, w, t)\mathcal{K}^S_{a,b}(\zeta, w)\right] dt\right| = \mathcal{O}(1). \tag{56}$$

Applying integration by parts and arguments from the proof of lemma 3.2 we represent the last integral as

$$\int_{\mathcal{U}(\epsilon)\times[0,1]} (1 - \phi_w(\zeta))\, D^c_w\left[c(\zeta, w, t)\mathcal{K}^S_{a,b}(\zeta, w)\right] dt \tag{57}$$

$$= \int_{\mathcal{U}(\epsilon)\times[0,1]} (1-\phi_w(\zeta))\left(D_w^c + D_\zeta^c\right)\left[c(\zeta,w,t)\mathcal{K}_{a,b}^S(\zeta,w)\right]dt$$

$$-\int_{\mathcal{U}(\epsilon)\times[0,1]} (1-\phi_w(\zeta))\, D_\zeta^c\left[c(\zeta,w,t)\mathcal{K}_{a,b}^S(\zeta,w)\right]dt$$

$$= \int_{\mathcal{U}(\epsilon)\times[0,1]} (1-\phi_w(\zeta))\left(D_w^c + D_\zeta^c\right)\left[c(\zeta,w,t)\mathcal{K}_{a,b}^S(\zeta,w)\right]dt$$

$$-\int_{b\mathcal{U}(\epsilon)\times[0,1]} c(\zeta,w,t)\cdot(1-\phi_w(\zeta))\cdot\mathcal{K}_{a,b}^S(\zeta,w)dt$$

$$+\sum_{I,d,j}\int_{\mathcal{U}(\epsilon)\times[0,1]} c_{\{I,d,j\}}(\zeta,w,t)\cdot(1-\phi_w(\zeta))\cdot\mathcal{K}_{d,h}^I(\zeta,w)dt$$

$$+\int_{\mathcal{U}(\epsilon)\times[0,1]} D_\zeta^c(1-\phi_w(\zeta))\left[c(\zeta,w,t)\mathcal{K}_{a,b}^S(\zeta,w)\right]dt$$

$$= -\int_{b\mathcal{U}(\epsilon)\times[0,1]} c(\zeta,w,t)\cdot(1-\phi_w(\zeta))\cdot\mathcal{K}_{a,b}^S(\zeta,w)dt$$

$$+\sum_{I,d,j}\int_{\mathcal{U}(\epsilon)\times[0,1]} c_{\{I,d,j\}}(\zeta,w,t)\cdot(1-\phi_w(\zeta))\cdot\mathcal{K}_{d,h}^I(\zeta,w)dt$$

$$+\int_{\mathcal{U}(\epsilon)\times[0,1]} D_\zeta^c(1-\phi_w(\zeta))\cdot c(\zeta,w,t)\mathcal{K}_{a,b}^S(\zeta,w)dt$$

$$+\sum_{\{L,e\}}\int_{\mathcal{U}(\epsilon)\times[0,1]} \left[Y_\zeta(w)(1-\phi_w(\zeta))\right]\cdot c_{\{L,e\}}(\zeta,w,t)\cdot\mathcal{K}_{e,b}^L(\zeta,w)dt,$$

with

$$\int_{b\mathcal{U}(\epsilon)\times[0,1]} c(\zeta,w,t)\cdot(1-\phi_w(\zeta))\cdot\mathcal{K}_{a,b}^S(\zeta,w)dt \in C^\infty(\mathcal{U}),$$

kernels $\mathcal{K}_{d,h}^I$ satisfying

$$k\left(\mathcal{K}_{d,h}^I\right) + h - l\left(\mathcal{K}_{d,h}^I\right) \leq k\left(\mathcal{K}_{a,b}^S\right) + b - l\left(\mathcal{K}_{a,b}^S\right),$$

$$k\left(\mathcal{K}_{d,h}^I\right) + 2h - 2l\left(\mathcal{K}_{d,h}^I\right) \leq k\left(\mathcal{K}_{a,b}^S\right) + 2b - 2l\left(\mathcal{K}_{a,b}^S\right),$$

and kernels $\mathcal{K}_{e,b}^L$ satisfying

$$k\left(\mathcal{K}_{e,b}^L\right) - l\left(\mathcal{K}_{e,b}^L\right) + 1 \leq k\left(\mathcal{K}_{a,b}^S\right) - l\left(\mathcal{K}_{a,b}^S\right),$$

$$k\left(\mathcal{K}_{e,b}^L\right) - 2l\left(\mathcal{K}_{e,b}^L\right) + 1 \leq k\left(\mathcal{K}_{a,b}^S\right) - 2l\left(\mathcal{K}_{a,b}^S\right). \tag{58}$$

For the second term of the right hand side of (57) using lemma 3.5 and conditions (42) we obtain

$$\left|\int_{\mathcal{U}(\epsilon)\times[0,1]} c_{\{I,d,j\}}(\zeta,w,t)\cdot(1-\phi_w(\zeta))\cdot\mathcal{K}_{d,h}^I(\zeta,w)dt\right|$$
$$= \mathcal{O}\left(\mathcal{I}_2\{0,k(\mathcal{K}),h(\mathcal{K})-l(\mathcal{K})\}(\epsilon,\delta)\right) = \mathcal{O}(1).$$

For the third term of the right hand side of (57) we obtain the following estimate

$$\left|\int_{\mathcal{U}(\epsilon)\times[0,1]} D_\zeta^c(1-\phi_w(\zeta))\cdot c(\zeta,w,t)\mathcal{K}_{a,b}^S(\zeta,w)dt\right|$$
$$= \delta^{-1}\cdot\mathcal{O}\left(\mathcal{I}_1\{0,k(\mathcal{K}),h(\mathcal{K})-l(\mathcal{K})\}(\epsilon,\delta)\right) = \mathcal{O}(1),$$

where we used properties of the function $\phi$, lemma 3.5 and conditions (42).

We obtain the same estimate for the integrals of the fourth term of the right hand side of (57) if we use the inequality

$$|\zeta - w| = \mathcal{O}(\delta) \text{ for } \zeta \in W'(w, \epsilon, \delta),$$

conditions (58) and estimate $|\mathrm{grad}\phi_w(\zeta)| = \mathcal{O}\left(1/\delta^2\right)$.

This completes the proof of the first estimate from (45) for the integral in (47). The only change that is necessary in this proof in order to make it work for the second estimate from (45) is the replacement $\delta \to \sqrt{\delta}$. Namely, instead of the neighborhoods $W(w, \epsilon, \delta)$ and $W'(w, \epsilon, \delta)$ we have to consider

$$\begin{aligned} W(w, \epsilon, \sqrt{\delta}) &= \left\{\zeta \in \widetilde{\mathcal{U}} : \rho(\zeta) = \epsilon, \quad |\Phi(\zeta, w)| \leq c\delta\right\}, \\ W'(w, \epsilon, \sqrt{\delta}) &= \left\{\zeta \in \widetilde{\mathcal{U}} : \rho(\zeta) = \epsilon, \quad |\Phi(\zeta, w)| \leq c'\delta\right\}, \end{aligned} \tag{59}$$

with $c < c'$ and such that if $|z - w| \leq \delta$ then

$$W'(w, \epsilon, \sqrt{\delta}) \subset W(z, \epsilon, C\sqrt{\delta})$$

with $c, c', C > 0$ independent of $z, w, \delta$.

Instead of function $\phi_w$ we have to consider a function $\psi_w \in C^\infty\left(\widetilde{\mathcal{U}} \setminus \mathcal{U}\right)$ such that $0 \leq \psi_w(\zeta) \leq 1$ and

$$\begin{aligned} \psi_w &\equiv 1 \text{ on } W(w, \epsilon, \sqrt{\delta}), \quad \psi_w \equiv 0 \text{ on } \mathcal{U}(\epsilon) \setminus W'(w, \epsilon, \sqrt{\delta}), \\ |D\psi_w(\zeta)| &= \mathcal{O}(1/\delta), \quad |D^c\psi_w(\zeta)| = \mathcal{O}(1/\sqrt{\delta}), \end{aligned} \tag{60}$$

for vector fields $D, D^c \in \mathbf{C}T(\mathbf{G})$ such that $D|_\mathcal{U} \in \mathbf{C}T(\mathbf{M})$ and $D^c|_\mathcal{U} \in \mathbf{C}T^c(\mathbf{M})$.

Accordingly, all the estimates should be changed with estimate (55) replaced by estimate (44). $\square$

Before proving the second estimate from (36) we will introduce additional notations.

Let $Y_1, \ldots, Y_m, U_1, \ldots, U_{n-m}, V_1, \ldots, V_{n-m}$ be tangent vector fields on $\widetilde{\mathcal{U}}$ such that for $\zeta \in \mathcal{U}$

(i) $\{Y_k(\zeta), (k = 1, \ldots, m), (U_j(\zeta), V_j(\zeta), (j = 1, \ldots, n - m)\}$ is a basis in $T_\zeta(\mathbf{M})$,

(ii) $\{(U_j(\zeta), V_j(\zeta), (j = 1, \ldots, n - m)\}$ is a basis in $T_\zeta^c(\mathbf{M})$.

For $z \in \mathcal{U}$ we define a diffeomorphism of a small neighborhood of the origin $\mathcal{V}_z \in T_z(\mathbf{G})$ onto the neighborhood $\widetilde{\mathcal{U}}$ by the formula

$$\sum_{k=1}^{m} x_k \frac{\partial}{\partial \rho_k}(z) + \sum_{k=1}^{m} y_k Y_k(z) + \sum_{j=1}^{n-m} u_j U_j(z) + \sum_{j=1}^{n-m} v_j V_j(z) \longrightarrow e_z(x, y, u, v)(1),$$

where $e_z(x, y, u, v)(t)$ is a solution of the system of differential equations

$$\frac{de_z}{dt}(x, y, u, v)(t) = \sum_{k=1}^{m} x_k \frac{\partial}{\partial \rho_k}(e_z(x, y, u, v)(t))$$

$$+ \sum_{k=1}^{m} y_k Y_k(e_z(x, y, u, v)(t)) + \sum_{j=1}^{n-m} u_j U_j(e_z(x, y, u, v)(t)) + \sum_{j=1}^{n-m} v_j V_j(e_z(x, y, u, v)(t))$$

with initial condition $e_z(x, y, u, v)(0) = z$.

For $z \in \mathcal{U}$ we also define a submanifold $\mathbf{M}_z^c \in \mathbf{M}$ as $e_z(\mathcal{V}_z \cap T_z^c(\mathbf{M}))$ and a map $\pi_z : \widetilde{\mathcal{U}} \to \mathbf{M}_z^c$ by the formula

$$\pi_z(\zeta) = e_z \circ p_z^c \circ e_z^{-1}(\zeta), \tag{61}$$

where $p_z^c$ is a projection of $T_z(\mathbf{G})$ onto $T_z^c(\mathbf{M})$ parallel to the vectors

$$\frac{\partial}{\partial \rho_1}(z), \ldots, \frac{\partial}{\partial \rho_m}(z), Y_1(z), \ldots, Y_m(z).$$

Then for any $\zeta \in \widetilde{\mathcal{U}}$ there exists a curve $x(\zeta, z) : [0, 1] \to \mathbf{M}_z^c$ defined by

$$x(\zeta, z, t) = e_z \left( p_z^c \circ e_z^{-1}(\zeta) \right)(t) \tag{62}$$

satisfying (2) with $x(\zeta, z, 0) = z, x(\zeta, z, 1) = \pi_z(\zeta)$.

In the lemma below we prove the second estimate from (36).

**Lemma 3.8.** *Let $0 < \alpha < 1$, $g \in \Pi^{1+\alpha}(\mathcal{U})$ be a function with compact support and $\widetilde{g} \in \Pi^{1+\alpha}\left(\{\rho\}, \widetilde{\mathcal{U}}(\epsilon_0)\right)$ its extension. Let $\mathcal{K}_{a,b}^S$ satisfy conditions (42).*

*Then*

$$f_\epsilon(z) := \left( \int_{\mathcal{U}(\epsilon) \times [0,1]} c(\zeta, z, t) \widetilde{g}(\zeta) \mathcal{K}_{a,b}^S(\zeta, z) dt \right) \in \Pi^{2+\alpha}(\mathcal{U})$$

*and*

$$\|f_\epsilon\|_{\Pi^{2+\alpha}(\mathcal{U})} \leq C \cdot \|g\|_{\Pi^{1+\alpha}(\mathcal{U})}$$

*with $C$ independent of $g$ and $\epsilon$.*

**Proof.** From the definition of spaces $\Pi^a$ we conclude that statement of the lemma would follow from inclusions

$$\begin{aligned} D^c f_\epsilon &\in \Lambda^{\frac{1+\alpha}{2}}(\mathcal{U}), \\ D^c D^c f_\epsilon, D f_\epsilon &\in \Gamma^\alpha(\mathcal{U}). \end{aligned} \tag{63}$$

We will start with the proof of the first of these inclusions. Using representation (46) we reduce the problem to the proof of inclusion

$$\int_{\mathcal{U}(\epsilon) \times [0,1]} (\widetilde{g}(\zeta) - g(z)) D_z^c \left[ c(\zeta, z, t) \mathcal{K}_{a,b}^S(\zeta, z) \right] dt \in \Lambda^{\frac{1+\alpha}{2}}(\mathcal{U}). \tag{64}$$

For a fixed $w \in \mathcal{U}$ and for $z$ such that $|z - w| \leq \delta$ we consider neighborhoods $W(w, \epsilon, \sqrt{\delta})$ and $W'(w, \epsilon, \sqrt{\delta})$ from (59), function $\psi_w(\zeta)$ from (60) and represent the integral from (64) as

$$\int_{\mathcal{U}(\epsilon) \times [0,1]} (\widetilde{g}(\zeta) - g(z)) D_z^c \left[ c(\zeta, z, t) \mathcal{K}_{a,b}^S(\zeta, z) \right] dt \tag{65}$$

$$= \int_{\mathcal{U}(\epsilon) \times [0,1]} (\widetilde{g}(\zeta) - g(\pi_z(\zeta))) \psi_w(\zeta) D_z^c \left[ c(\zeta, z, t) \mathcal{K}_{a,b}^S(\zeta, z) \right] dt$$

$$+ \int_{\mathcal{U}(\epsilon) \times [0,1]} (\widetilde{g}(\zeta) - g(\pi_z(\zeta))) (1 - \psi_w(\zeta)) D_z^c \left[ c(\zeta, z, t) \mathcal{K}_{a,b}^S(\zeta, z) \right] dt$$

$$+ \int_{\mathcal{U}(\epsilon) \times [0,1]} (g(\pi_z(\zeta)) - g(z)) D_z^c \left[ c(\zeta, z, t) \mathcal{K}_{a,b}^S(\zeta, z) \right] dt,$$

with $\pi_z(\zeta)$ defined in (61).

Then for the first integral in the right hand side of (65) using estimates

$$|\widetilde{g}(\zeta) - g(\pi_z(\zeta))| = \|g\|_{\Gamma^{1+\alpha}(\mathcal{U})} \cdot \mathcal{O}\left(|\Phi(\zeta, z)|^{\frac{1+\alpha}{2}}\right) \tag{66}$$

and
$$|\Phi(\zeta, z)| = \mathcal{O}(\delta)$$
for $\zeta \in W(z, \epsilon, C\sqrt{\delta})$, representation (52) with kernels satisfying (53) and lemma 3.5 we obtain the estimate

$$\left| \int_{\mathcal{U}(\epsilon) \times [0,1]} (\widetilde{g}(\zeta) - g(\pi_z(\zeta))) \psi_w(\zeta) c(\zeta, z, t) \mathcal{K}^I_{d,h}(\zeta, z) dt \right|$$

$$= \mathcal{O}\left( \left| \int_{W(z,\epsilon,C\sqrt{\delta}) \times [0,1]} (\widetilde{g}(\zeta) - g(\pi_z(\zeta))) c(\zeta, z, t) \mathcal{K}^I_{d,h}(\zeta, z) dt \right| \right)$$

$$= \|g\|_{\Gamma^{1+\alpha}(\mathcal{U})} \cdot \delta^{\frac{\alpha}{2}} \cdot \mathcal{O}\left( \int_{W(z,\epsilon,C\sqrt{\delta}) \times [0,1]} \frac{\overbrace{\wedge_i d\theta_i(\zeta)}^{m-1} \wedge d\sigma_{2n-m}(\zeta)}{|\zeta - z|^{k(\mathcal{K})} \cdot |\Phi(\zeta, z)|^{h(\mathcal{K}) - l(\mathcal{K}) - 1/2}} \right)$$

$$= \|g\|_{\Gamma^{1+\alpha}(\mathcal{U})} \cdot \delta^{\frac{\alpha}{2}} \cdot \mathcal{O}\left( \mathcal{I}_1 \left\{ 0, k(\mathcal{K}), h(\mathcal{K}) - l(\mathcal{K}) - \frac{1}{2} \right\} (\epsilon, \sqrt{\delta}) \right)$$

$$= \|g\|_{\Gamma^{1+\alpha}(\mathcal{U})} \cdot \mathcal{O}\left( \delta^{\frac{1+\alpha}{2}} \right).$$

For the second integral in the right hand side of (65) we consider representation

$$\int_{\mathcal{U}(\epsilon) \times [0,1]} (\widetilde{g}(\zeta) - g(\pi_z(\zeta))) (1 - \psi_w(\zeta)) D^c_z \left[ c(\zeta, z, t) \mathcal{K}^S_{a,b}(\zeta, z) \right] dt \qquad (67)$$

$$- \int_{\mathcal{U}(\epsilon) \times [0,1]} (\widetilde{g}(\zeta) - g(\pi_w(\zeta))) (1 - \psi_w(\zeta)) D^c_w \left[ c(\zeta, w, t) \mathcal{K}^S_{a,b}(\zeta, w) \right] dt$$

$$= \int_{\mathcal{U}(\epsilon) \times [0,1]} (\widetilde{g}(\zeta) - g(\pi_z(\zeta))) (1 - \psi_w(\zeta))$$

$$\times \left( D^c_z \left[ c(\zeta, z, t) \mathcal{K}^S_{a,b}(\zeta, z) \right] - D^c_w \left[ c(\zeta, w, t) \mathcal{K}^S_{a,b}(\zeta, w) \right] \right) dt$$

$$+ \int_{\mathcal{U}(\epsilon) \times [0,1]} (g(\pi_w(\zeta)) - g(\pi_z(\zeta))) (1 - \psi_w(\zeta)) D^c_w \left[ c(\zeta, w, t) \mathcal{K}^S_{a,b}(\zeta, w) \right] dt.$$

Then for the first term of the right hand side of (67) using representation (52) with kernels satisfying (53), estimates (44) and (66) and lemma 3.5 we obtain

$$\left| \int_{\mathcal{U}(\epsilon) \times [0,1]} (\widetilde{g}(\zeta) - g(\pi_z(\zeta))) (1 - \psi_w(\zeta)) \right.$$

$$\left. \times \left( D^c_z \left[ c(\zeta, z, t) \mathcal{K}^S_{a,b}(\zeta, z) \right] - D^c_w \left[ c(\zeta, w, t) \mathcal{K}^S_{a,b}(\zeta, w) \right] \right) dt \right|$$

$$= \left| \sum_{I,d,j} \int_{\mathcal{U}(\epsilon) \times [0,1]} (\widetilde{g}(\zeta) - g(\pi_z(\zeta))) (1 - \psi_w(\zeta)) \right.$$

$$\left. \times \left[ c_{\{I,d,j\}}(\zeta, z, t) \cdot \mathcal{K}^I_{d,h}(\zeta, z) - c_{\{I,d,j\}}(\zeta, w, t) \cdot \mathcal{K}^I_{d,h}(\zeta, w) \right] dt \right|$$

$$= \|g\|_{\Gamma^{1+\alpha}(\mathcal{U})} \cdot \left[ \delta \cdot \mathcal{O}\left( \mathcal{I}_2 \{\alpha, k(\mathcal{K}), h(\mathcal{K}) - l(\mathcal{K}) - 1/2\} (\epsilon, \sqrt{\delta}) \right) \right.$$

$$+ \delta \cdot \mathcal{O}\left( \mathcal{I}_2 \{\alpha, k(\mathcal{K}) + 1, h(\mathcal{K}) - l(\mathcal{K}) - 1/2\} (\epsilon, \sqrt{\delta}) \right)$$

$$\left. + \delta \cdot \mathcal{O}\left( \mathcal{I}_2 \{\alpha, k(\mathcal{K}), h(\mathcal{K}) - l(\mathcal{K}) + 1/2\} (\epsilon, \sqrt{\delta}) \right) \right] = \|g\|_{\Gamma^{1+\alpha}(\mathcal{U})} \cdot \mathcal{O}\left( \delta^{\frac{1+\alpha}{2}} \right).$$

For the second term of the right hand side of (67) we use representation

$$\int_{\mathcal{U}(\epsilon) \times [0,1]} (g(\pi_w(\zeta)) - g(\pi_z(\zeta))) (1 - \psi_w(\zeta)) D^c_w \left[ c(\zeta, w, t) \mathcal{K}^S_{a,b}(\zeta, w) \right] dt \qquad (68)$$

$$= (g(w) - g(z)) \cdot \int_{\mathcal{U}(\epsilon) \times [0,1]} (1 - \psi_w(\zeta)) \, D_w^c \left[ c(\zeta, w, t) \mathcal{K}_{a,b}^S(\zeta, w) \right] dt$$

$$+ \int_{\mathcal{U}(\epsilon) \times [0,1]} (g(\pi_w(\zeta)) - g(w)) \, (1 - \psi_w(\zeta)) \, D_w^c \left[ c(\zeta, w, t) \mathcal{K}_{a,b}^S(\zeta, w) \right] dt$$

$$- \int_{\mathcal{U}(\epsilon) \times [0,1]} (g(\pi_z(\zeta)) - g(z)) \, (1 - \psi_w(\zeta)) \, D_w^c \left[ c(\zeta, w, t) \mathcal{K}_{a,b}^S(\zeta, w) \right] dt,$$

with the first integral in the right hand side of (68) admitting the estimate

$$\|g\|_{\Gamma^{1+\alpha}(\mathcal{U})} \cdot \mathcal{O}\left(|w - z|^{\frac{1+\alpha}{2}}\right)$$

which follows from the estimate

$$|g(w) - g(z)| = \|g\|_{\Gamma^{1+\alpha}(\mathcal{U})} \cdot \mathcal{O}\left(|w - z|^{\frac{1+\alpha}{2}}\right)$$

and the estimate

$$\left| \int_{\mathcal{U}(\epsilon) \times [0,1]} (1 - \psi_w(\zeta)) \, D_w^c \left[ c(\zeta, w, t) \mathcal{K}_{a,b}^S(\zeta, w) \right] dt \right| = \mathcal{O}(1),$$

analogous to the estimate (56).

Joining now the last two integrals from the right hand side of (68) with the last integral from (65) we reduce the problem to the estimate of the integral

$$\int_{\mathcal{U}(\epsilon) \times [0,1]} (g(\pi_w(\zeta)) - g(w)) \, (1 - \psi_w(\zeta)) \, D_w^c \left[ c(\zeta, w, t) \mathcal{K}_{a,b}^S(\zeta, w) \right] dt \qquad (69)$$

$$- \int_{\mathcal{U}(\epsilon) \times [0,1]} (g(\pi_z(\zeta)) - g(z)) \, (1 - \psi_w(\zeta)) \, D_w^c \left[ c(\zeta, w, t) \mathcal{K}_{a,b}^S(\zeta, w) \right] dt$$

$$+ \int_{\mathcal{U}(\epsilon) \times [0,1]} (g(\pi_z(\zeta)) - g(z)) \, D_z^c \left[ c(\zeta, z, t) \mathcal{K}_{a,b}^S(\zeta, z) \right] dt$$

$$- \int_{\mathcal{U}(\epsilon) \times [0,1]} (g(\pi_w(\zeta)) - g(w)) \, D_w^c \left[ c(\zeta, w, t) \mathcal{K}_{a,b}^S(\zeta, w) \right] dt$$

$$= \int_{\mathcal{U}(\epsilon) \times [0,1]} (g(\pi_z(\zeta)) - g(z)) \, (1 - \psi_w(\zeta))$$

$$\times \left( D_z^c \left[ c(\zeta, z, t) \mathcal{K}_{a,b}^S(\zeta, z) \right] - D_w^c \left[ c(\zeta, w, t) \mathcal{K}_{a,b}^S(\zeta, w) \right] \right) dt$$

$$+ \int_{\mathcal{U}(\epsilon) \times [0,1]} (g(\pi_z(\zeta)) - g(z)) \cdot \psi_w(\zeta) D_z^c \left[ c(\zeta, z, t) \mathcal{K}_{a,b}^S(\zeta, z) \right] dt$$

$$- \int_{\mathcal{U}(\epsilon) \times [0,1]} (g(\pi_w(\zeta)) - g(w)) \cdot \psi_w(\zeta) D_w^c \left[ c(\zeta, w, t) \mathcal{K}_{a,b}^S(\zeta, w) \right] dt.$$

Integrating along the curve $x(\zeta, z, s)$ from $z$ to $\pi_z(\zeta)$ and along the curve $x(\zeta, w, s)$ from $w$ to $\pi_w(\zeta)$ we transform the right hand side of (69):

$$\int_{\mathcal{U}(\epsilon) \times [0,1]} (g(\pi_z(\zeta)) - g(z)) \, (1 - \psi_w(\zeta)) \qquad (70)$$

$$\times \left( D_z^c \left[ c(\zeta, z, t) \mathcal{K}_{a,b}^S(\zeta, z) \right] - D_w^c \left[ c(\zeta, w, t) \mathcal{K}_{a,b}^S(\zeta, w) \right] \right) dt$$

$$+ \int_{\mathcal{U}(\epsilon) \times [0,1]} (g(\pi_z(\zeta)) - g(z)) \cdot \psi_w(\zeta) D_z^c \left[ c(\zeta, z, t) \mathcal{K}_{a,b}^S(\zeta, z) \right] dt$$

$$- \int_{\mathcal{U}(\epsilon) \times [0,1]} (g(\pi_w(\zeta)) - g(w)) \cdot \psi_w(\zeta) D_w^c \left[ c(\zeta, w, t) \mathcal{K}_{a,b}^S(\zeta, w) \right] dt$$

$$= \int_{[0,1]} ds \int_{\mathcal{U}(\epsilon) \times [0,1]} \langle \mathrm{grad}_c g(x(\zeta, z, s)), x'(\zeta, z, s) \rangle \, (1 - \psi_w(\zeta))$$

$$\times \left( D_z^c \left[ c(\zeta,z,t) \mathcal{K}_{a,b}^S(\zeta,z) \right] - D_w^c \left[ c(\zeta,w,t) \mathcal{K}_{a,b}^S(\zeta,w) \right] \right) dt$$

$$+ \int_{[0,1]} ds \int_{\mathcal{U}(\epsilon) \times [0,1]} \langle \mathrm{grad}_c g(x(\zeta,z,s)), x'(\zeta,z,s) \rangle \cdot \psi_w(\zeta) D_z^c \left[ c(\zeta,z,t) \mathcal{K}_{a,b}^S(\zeta,z) \right] dt$$

$$- \int_{[0,1]} ds \int_{\mathcal{U}(\epsilon) \times [0,1]} \langle \mathrm{grad}_c g(x(\zeta,w,s)), x'(\zeta,w,s) \rangle \cdot \psi_w(\zeta) D_w^c \left[ c(\zeta,w,t) \mathcal{K}_{a,b}^S(\zeta,w) \right] dt$$

$$= \int_{[0,1]} ds \int_{\mathcal{U}(\epsilon) \times [0,1]} \langle \mathrm{grad}_c g(x(\zeta,z,s)) - \mathrm{grad}_c g(z), x'(\zeta,z,s) \rangle (1 - \psi_w(\zeta))$$

$$\times \left( D_z^c \left[ c(\zeta,z,t) \mathcal{K}_{a,b}^S(\zeta,z) \right] - D_w^c \left[ c(\zeta,w,t) \mathcal{K}_{a,b}^S(\zeta,w) \right] \right) dt$$

$$+ \int_{[0,1]} ds \int_{\mathcal{U}(\epsilon) \times [0,1]} \langle \mathrm{grad}_c g(x(\zeta,z,s)) - \mathrm{grad}_c g(z), x'(\zeta,z,s) \rangle$$

$$\times \psi_w(\zeta) D_z^c \left[ c(\zeta,z,t) \mathcal{K}_{a,b}^S(\zeta,z) \right] dt$$

$$- \int_{[0,1]} ds \int_{\mathcal{U}(\epsilon) \times [0,1]} \langle \mathrm{grad}_c g(x(\zeta,w,s)) - \mathrm{grad}_c g(w), x'(\zeta,w,s) \rangle$$

$$\times \psi_w(\zeta) D_w^c \left[ c(\zeta,w,t) \mathcal{K}_{a,b}^S(\zeta,w) \right] dt$$

$$+ \int_{[0,1]} ds \int_{\mathcal{U}(\epsilon) \times [0,1]} \langle \mathrm{grad}_c g(z), x'(\zeta,z,s) \rangle (1 - \psi_w(\zeta))$$

$$\times \left( D_z^c \left[ c(\zeta,z,t) \mathcal{K}_{a,b}^S(\zeta,z) \right] - D_w^c \left[ c(\zeta,w,t) \mathcal{K}_{a,b}^S(\zeta,w) \right] \right) dt$$

$$+ \int_{[0,1]} ds \int_{\mathcal{U}(\epsilon) \times [0,1]} \langle \mathrm{grad}_c g(z), x'(\zeta,z,s) \rangle \cdot \psi_w(\zeta) D_z^c \left[ c(\zeta,z,t) \mathcal{K}_{a,b}^S(\zeta,z) \right] dt$$

$$- \int_{[0,1]} ds \int_{\mathcal{U}(\epsilon) \times [0,1]} \langle \mathrm{grad}_c g(w), x'(\zeta,w,s) \rangle \cdot \psi_w(\zeta) D_w^c \left[ c(\zeta,w,t) \mathcal{K}_{a,b}^S(\zeta,w) \right] dt,$$

where $\mathrm{grad}_c g$ is the projection of the gradient of $g$ on $T^c(\mathbf{M})$ and $x'(\zeta,z,s)$ is a vector with components $dx_i(\zeta,z,s)/ds$.

Indices of kernels $(dx_i(\zeta,z,s)/ds) \cdot D_z^c \left[ c(\zeta,z,t) \mathcal{K}_{a,b}^S(\zeta,z) \right]$ satisfy conditions (42) and

$$|\mathrm{grad}_c g(x(\zeta,z,s)) - \mathrm{grad}_c g(z)| = \|g\|_{\Pi^{1+\alpha}(\mathcal{U})} \cdot \mathcal{O}\left( |\pi_z(\zeta) - z|^\alpha \right) \tag{71}$$

$$= \|g\|_{\Pi^{1+\alpha}(\mathcal{U})} \cdot \mathcal{O}\left( |\Phi(\zeta,z)|^{\frac{\alpha}{2}} \right).$$

Then for the first term of the right hand side of (70) using estimates (44), (71) and lemma 3.5 we obtain

$$\left| \int_{[0,1]} ds \int_{\mathcal{U}(\epsilon) \times [0,1]} \langle \mathrm{grad}_c g(x(\zeta,z,s)) - \mathrm{grad}_c g(z), x'(\zeta,z,s) \rangle (1 - \psi_w(\zeta)) \right.$$

$$\left. \times \left( D_z^c \left[ c(\zeta,z,t) \mathcal{K}_{a,b}^S(\zeta,z) \right] - D_w^c \left[ c(\zeta,w,t) \mathcal{K}_{a,b}^S(\zeta,w) \right] \right) dt \right|$$

$$= \|g\|_{\Gamma^{1+\alpha}(\mathcal{U})} \cdot \left[ \delta \cdot \mathcal{O}\left( \mathcal{I}_2 \{\alpha, k(\mathcal{K}), h(\mathcal{K}) - l(\mathcal{K})\} \left(\epsilon, \sqrt{\delta}\right) \right) \right.$$

$$+ \delta \cdot \mathcal{O}\left( \mathcal{I}_2 \{\alpha, k(\mathcal{K}) + 1, h(\mathcal{K}) - l(\mathcal{K})\} \left(\epsilon, \sqrt{\delta}\right) \right)$$

$$\left. + \delta \cdot \mathcal{O}\left( \mathcal{I}_2 \{\alpha, k(\mathcal{K}), h(\mathcal{K}) - l(\mathcal{K}) + 1\} \left(\epsilon, \sqrt{\delta}\right) \right) \right]$$

$$= \|g\|_{\Gamma^{1+\alpha}(\mathcal{U})} \cdot \delta \cdot \mathcal{O}\left( \delta^{\frac{\alpha-1}{2}} \right) = \|g\|_{\Gamma^{1+\alpha}(\mathcal{U})} \cdot \mathcal{O}\left( \delta^{\frac{1+\alpha}{2}} \right).$$

Again using the fact that the indices of kernels $(dx_i(\zeta,z,s)/ds) \cdot D_z^c \left[ c(\zeta,z,t)\mathcal{K}_{a,b}^S(\zeta,z) \right]$ satisfy conditions (42), estimate (71) and lemma 3.5 we obtain for the second and third terms of the right hand side of (70)

$$\left| \int_{[0,1]} ds \int_{\mathcal{U}(\epsilon)\times[0,1]} \langle \operatorname{grad}_c g(x(\zeta,z,s)) - \operatorname{grad}_c g(z), x'(\zeta,z,s) \rangle \right.$$
$$\left. \times \psi_w(\zeta) D_z^c \left[ c(\zeta,z,t)\mathcal{K}_{a,b}^S(\zeta,z) \right] dt \right|$$
$$= \|g\|_{\Gamma^{1+\alpha}(\mathcal{U})} \cdot \delta^{\frac{\alpha}{2}} \cdot \mathcal{O} \left( \int_{W(z,\epsilon,C\sqrt{\delta})\times[0,1]} \left| D_z^c \left[ c(\zeta,z,t)\mathcal{K}_{a,b}^S(\zeta,z) \right] \right| dt \right)$$
$$= \|g\|_{\Gamma^{1+\alpha}(\mathcal{U})} \cdot \delta^{\frac{\alpha}{2}} \cdot \mathcal{O} \left( \mathcal{I}_1 \{0, k(\mathcal{K}), h(\mathcal{K}) - l(\mathcal{K})\} \left( \epsilon, \sqrt{\delta} \right) \right)$$
$$= \|g\|_{\Gamma^{1+\alpha}(\mathcal{U})} \cdot \mathcal{O} \left( \delta^{\frac{1+\alpha}{2}} \right).$$

For the rest of the integrals from the right hand side of (70) we use representation

$$\int_{[0,1]} ds \int_{\mathcal{U}(\epsilon)\times[0,1]} \langle \operatorname{grad}_c g(z), x'(\zeta,z,s) \rangle (1 - \psi_w(\zeta)) \quad (72)$$
$$\times \left( D_z^c \left[ c(\zeta,z,t)\mathcal{K}_{a,b}^S(\zeta,z) \right] - D_w^c \left[ c(\zeta,w,t)\mathcal{K}_{a,b}^S(\zeta,w) \right] \right) dt$$
$$+ \int_{[0,1]} ds \int_{\mathcal{U}(\epsilon)\times[0,1]} \langle \operatorname{grad}_c g(z), x'(\zeta,z,s) \rangle \cdot \psi_w(\zeta) D_z^c \left[ c(\zeta,z,t)\mathcal{K}_{a,b}^S(\zeta,z) \right] dt$$
$$- \int_{[0,1]} ds \int_{\mathcal{U}(\epsilon)\times[0,1]} \langle \operatorname{grad}_c g(w), x'(\zeta,w,s) \rangle \cdot \psi_w(\zeta) D_w^c \left[ c(\zeta,w,t)\mathcal{K}_{a,b}^S(\zeta,w) \right] dt$$
$$= \int_{[0,1]} ds \int_{\mathcal{U}(\epsilon)\times[0,1]} \langle \operatorname{grad}_c g(z) - \operatorname{grad}_c g(w), x'(\zeta,z,s) \rangle$$
$$\times \psi_w(\zeta) D_z^c \left[ c(\zeta,z,t)\mathcal{K}_{a,b}^S(\zeta,z) \right] dt$$
$$+ \left( \int_{[0,1]} ds \int_{\mathcal{U}(\epsilon)\times[0,1]} \langle \operatorname{grad}_c g(w), x'(\zeta,z,s) \rangle D_z^c \left[ c(\zeta,z,t)\mathcal{K}_{a,b}^S(\zeta,z) \right] dt \right.$$
$$\left. - \int_{[0,1]} ds \int_{\mathcal{U}(\epsilon)\times[0,1]} \langle \operatorname{grad}_c g(w), x'(\zeta,w,s) \rangle D_w^c \left[ c(\zeta,w,t)\mathcal{K}_{a,b}^S(\zeta,w) \right] dt \right)$$
$$- \int_{[0,1]} ds \int_{\mathcal{U}(\epsilon)\times[0,1]} \langle \operatorname{grad}_c g(w), x'(\zeta,z,s) - x'(\zeta,w,s) \rangle$$
$$\times (1 - \psi_w(\zeta)) D_z^c \left[ c(\zeta,z,t)\mathcal{K}_{a,b}^S(\zeta,z) \right] dt$$
$$+ \int_{[0,1]} ds \int_{\mathcal{U}(\epsilon)\times[0,1]} \langle \operatorname{grad}_c g(z), x'(\zeta,z,s) \rangle (1 - \psi_w(\zeta))$$
$$\times \left( D_z^c \left[ c(\zeta,z,t)\mathcal{K}_{a,b}^S(\zeta,z) \right] - D_w^c \left[ c(\zeta,w,t)\mathcal{K}_{a,b}^S(\zeta,w) \right] \right) dt$$
$$- \int_{[0,1]} ds \int_{\mathcal{U}(\epsilon)\times[0,1]} \langle \operatorname{grad}_c g(w), x'(\zeta,w,s) \rangle (1 - \psi_w(\zeta))$$
$$\times \left( D_z^c \left[ c(\zeta,z,t)\mathcal{K}_{a,b}^S(\zeta,z) \right] - D_w^c \left[ c(\zeta,w,t)\mathcal{K}_{a,b}^S(\zeta,w) \right] \right) dt$$

Using estimate
$$|\operatorname{grad}_c g(z) - \operatorname{grad}_c g(w)| = \|g\|_{\Pi^{1+\alpha}(\mathcal{U})} \cdot \mathcal{O}\left(\delta^{\frac{\alpha}{2}}\right), \quad (73)$$

and the fact that indices of kernels $(dx_i(\zeta,z,s)/ds) \cdot D_z^c \left[ c(\zeta,z,t)\mathcal{K}_{a,b}^S(\zeta,z) \right]$ satisfy conditions (42) we obtain for the first integral in the right hand side of (72) an estimate $\|g\|_{\Pi^{1+\alpha}(\mathcal{U})} \cdot$

$\mathcal{O}\left(\delta^{\frac{1+\alpha}{2}}\right).$

According to the arguments in the first part of lemma 3.7

$$\int_{[0,1]} ds \int_{\mathcal{U}(\epsilon)\times[0,1]} \langle \operatorname{grad}_c g(w), x'(\zeta,z,s)\rangle D_z^c\left[c(\zeta,z,t)\mathcal{K}_{a,b}^S(\zeta,z)\right] dt$$

is a function from $\Lambda^1(\mathcal{U})$ therefore for the second term in the right hand side of (72) we have

$$\left|\int_{[0,1]} ds \int_{\mathcal{U}(\epsilon)\times[0,1]} \langle \operatorname{grad}_c g(w), x'(\zeta,z,s)\rangle D_z^c\left[c(\zeta,z,t)\mathcal{K}_{a,b}^S(\zeta,z)\right] dt \right.$$
$$\left. - \int_{[0,1]} ds \int_{\mathcal{U}(\epsilon)\times[0,1]} \langle \operatorname{grad}_c g(w), x'(\zeta,w,s)\rangle D_w^c\left[c(\zeta,w,t)\mathcal{K}_{a,b}^S(\zeta,w)\right] dt \right| = \mathcal{O}(\delta).$$

Using estimate
$$|x'(\zeta,z,s) - x'(\zeta,w,s)| = \mathcal{O}(\delta), \tag{74}$$

representation (52) for the kernel $D_z^c\left[c(\zeta,z,t)\mathcal{K}_{a,b}^S(\zeta,z)\right]$ with kernels satisfying conditions (53) and lemma 3.5 we obtain for the third term of the right side of (72)

$$\left|\int_{\mathcal{U}(\epsilon)\times[0,1]} \langle \operatorname{grad}_c g(w), x'(\zeta,z,s) - x'(\zeta,w,s)\rangle (1-\psi_w(\zeta)) D_z^c\left[c(\zeta,z,t)\mathcal{K}_{a,b}^S(\zeta,z)\right] dt\right|$$

$$= \|g\|_{\Gamma^{1+\alpha}(\mathcal{U})} \cdot \delta \cdot \mathcal{O}\left(\int_{\mathcal{U}(\epsilon)\setminus W_w'(\epsilon,\sqrt\delta)} \frac{\overbrace{\wedge_i d\theta_i(\zeta)}^{m-1} \wedge d\sigma_{2n-m}(\zeta)}{|\zeta-z|^{k(\mathcal{K})} \cdot |\Phi(\zeta,z)|^{h(\mathcal{K})-l(\mathcal{K})}}\right)$$

$$= \|g\|_{\Gamma^{1+\alpha}(\mathcal{U})} \cdot \delta \cdot \mathcal{O}\left(\mathcal{I}_2\{0,k(\mathcal{K}),h(\mathcal{K})-l(\mathcal{K})\}\left(\epsilon,\sqrt\delta\right)\right)$$
$$= \|g\|_{\Gamma^{1+\alpha}(\mathcal{U})} \cdot \mathcal{O}(\delta\log\delta).$$

For the last two integrals from the right hand side of (72) we use representation

$$\int_{[0,1]} ds \int_{\mathcal{U}(\epsilon)\times[0,1]} \langle \operatorname{grad}_c g(z), x'(\zeta,z,s)\rangle (1-\psi_w(\zeta)) \tag{75}$$
$$\times \left(D_z^c\left[c(\zeta,z,t)\mathcal{K}_{a,b}^S(\zeta,z)\right] - D_w^c\left[c(\zeta,w,t)\mathcal{K}_{a,b}^S(\zeta,w)\right]\right) dt$$
$$- \int_{[0,1]} ds \int_{\mathcal{U}(\epsilon)\times[0,1]} \langle \operatorname{grad}_c g(w), x'(\zeta,w,s)\rangle (1-\psi_w(\zeta))$$
$$\times \left(D_z^c\left[c(\zeta,z,t)\mathcal{K}_{a,b}^S(\zeta,z)\right] - D_w^c\left[c(\zeta,w,t)\mathcal{K}_{a,b}^S(\zeta,w)\right]\right) dt$$
$$= \int_{[0,1]} ds \int_{\mathcal{U}(\epsilon)\times[0,1]} \langle \operatorname{grad}_c g(z) - \operatorname{grad}_c g(w), x'(\zeta,z,s)\rangle (1-\psi_w(\zeta))$$
$$\times \left(D_z^c\left[c(\zeta,z,t)\mathcal{K}_{a,b}^S(\zeta,z)\right] - D_w^c\left[c(\zeta,w,t)\mathcal{K}_{a,b}^S(\zeta,w)\right]\right) dt$$
$$+ \int_{[0,1]} ds \int_{\mathcal{U}(\epsilon)\times[0,1]} \langle \operatorname{grad}_c g(w), x'(\zeta,z,s) - x'(\zeta,w,s)\rangle (1-\psi_w(\zeta))$$
$$\times \left(D_z^c\left[c(\zeta,z,t)\mathcal{K}_{a,b}^S(\zeta,z)\right] - D_w^c\left[c(\zeta,w,t)\mathcal{K}_{a,b}^S(\zeta,w)\right]\right) dt.$$

Then for the first term of the right hand side of (75) using representation (52) with kernels satisfying (53), estimates (44), (73) and lemma 3.5 we obtain

$$\left|\int_{[0,1]} ds \int_{\mathcal{U}(\epsilon)\times[0,1]} \langle \operatorname{grad}_c g(z) - \operatorname{grad}_c g(w), x'(\zeta,z,s)\rangle (1-\psi_w(\zeta)) \right.$$
$$\left. \times \left(D_z^c\left[c(\zeta,z,t)\mathcal{K}_{a,b}^S(\zeta,z)\right] - D_w^c\left[c(\zeta,w,t)\mathcal{K}_{a,b}^S(\zeta,w)\right]\right) dt\right|$$

$$= \|g\|_{\Gamma^{1+\alpha}(\mathcal{U})} \cdot \delta^{\frac{\alpha}{2}} \cdot \left[ \delta \cdot \mathcal{O}\left(\mathcal{I}_2\{0, k(\mathcal{K})+1, h(\mathcal{K})-l(\mathcal{K})\}\left(\epsilon, \sqrt{\delta}\right)\right) \right.$$
$$+\delta \cdot \mathcal{O}\left(\mathcal{I}_2\{0, k(\mathcal{K})-1, h(\mathcal{K})-l(\mathcal{K})+1\}\left(\epsilon, \sqrt{\delta}\right)\right)$$
$$+\delta \cdot \mathcal{O}\left(\mathcal{I}_2\{0, k(\mathcal{K}), h(\mathcal{K})-l(\mathcal{K})+1\}\left(\epsilon, \sqrt{\delta}\right)\right)$$
$$\left.+\delta \cdot \mathcal{O}\left(\mathcal{I}_2\{0, k(\mathcal{K})-2, h(\mathcal{K})-l(\mathcal{K})+2\}\left(\epsilon, \sqrt{\delta}\right)\right)\right]$$
$$= \|g\|_{\Gamma^{1+\alpha}(\mathcal{U})} \cdot \delta^{\frac{\alpha}{2}+1} \cdot \mathcal{O}\left(\delta^{-\frac{1}{2}}\right) = \|g\|_{\Gamma^{1+\alpha}(\mathcal{U})} \cdot \mathcal{O}\left(\delta^{\frac{1+\alpha}{2}}\right).$$

For the second term of the right hand side of (75) using representation (52) with kernels satisfying (53), estimates (44), (74) and lemma 3.5 we obtain

$$\left| \int_{[0,1]} ds \int_{\mathcal{U}(\epsilon)\times[0,1]} \langle \mathrm{grad}_c g(w), x'(\zeta,z,s) - x'(\zeta,w,s) \rangle (1-\psi_w(\zeta)) \right.$$
$$\left. \times \left( D_z^c \left[ c(\zeta,z,t) \mathcal{K}_{a,b}^S(\zeta,z) \right] - D_w^c \left[ c(\zeta,w,t) \mathcal{K}_{a,b}^S(\zeta,w) \right] \right) dt \right|$$
$$= \|g\|_{\Gamma^{1+\alpha}(\mathcal{U})} \cdot \delta \cdot \left[ \delta \cdot \mathcal{O}\left(\mathcal{I}_2\{0, k(\mathcal{K})+2, h(\mathcal{K})-l(\mathcal{K})\}\left(\epsilon, \sqrt{\delta}\right)\right) \right.$$
$$+\delta \cdot \mathcal{O}\left(\mathcal{I}_2\{0, k(\mathcal{K}), h(\mathcal{K})-l(\mathcal{K})+1\}\left(\epsilon, \sqrt{\delta}\right)\right)$$
$$+\delta \cdot \mathcal{O}\left(\mathcal{I}_2\{0, k(\mathcal{K})+1, h(\mathcal{K})-l(\mathcal{K})+1\}\left(\epsilon, \sqrt{\delta}\right)\right)$$
$$\left.+\delta \cdot \mathcal{O}\left(\mathcal{I}_2\{0, k(\mathcal{K})-1, h(\mathcal{K})-l(\mathcal{K})+2\}\left(\epsilon, \sqrt{\delta}\right)\right)\right]$$
$$= \|g\|_{\Gamma^{1+\alpha}(\mathcal{U})} \cdot \delta^2 \cdot \mathcal{O}\left(\delta^{-1}\right) = \|g\|_{\Gamma^{1+\alpha}(\mathcal{U})} \cdot \mathcal{O}(\delta).$$

To prove inclusion $D^c D^c f_\epsilon \in \Lambda^{\frac{\alpha}{2}}(\mathcal{U})$ we consider representation

$$D^c D^c f_\epsilon(z) = g(z) \cdot D_z^c D_z^c \int_{\mathcal{U}(\epsilon)\times[0,1]} c(\zeta,z,t) \mathcal{K}_{a,b}^S(\zeta,z) dt \qquad (76)$$
$$+ \langle \mathrm{grad}_c g(z), D_z^c D_z^c \int_{\mathcal{U}(\epsilon)\times[0,1]} (\zeta-z) \cdot c(\zeta,z,t) \mathcal{K}_{a,b}^S(\zeta,z) dt \rangle$$
$$- 2\langle \mathrm{grad}_c g(z), D_z^c \int_{\mathcal{U}(\epsilon)\times[0,1]} [D_z^c(\zeta-z)] \cdot c(\zeta,z,t) \mathcal{K}_{a,b}^S(\zeta,z) dt \rangle$$
$$+ \langle \mathrm{grad}_c g(z), \int_{\mathcal{U}(\epsilon)\times[0,1]} [D_z^c D_z^c(\zeta-z)] \cdot c(\zeta,z,t) \mathcal{K}_{a,b}^S(\zeta,z) dt \rangle$$
$$+ \int_{\mathcal{U}(\epsilon)\times[0,1]} [\widetilde{g}(\zeta) - g(z) - \langle \mathrm{grad}_c g(z), \zeta-z \rangle] D_z^c D_z^c \left[ c(\zeta,z,t) \mathcal{K}_{a,b}^S(\zeta,z) \right] dt$$

and reduce the problem to the proof of inclusion

$$\int_{\mathcal{U}(\epsilon)\times[0,1]} [\widetilde{g}(\zeta) - g(z) - \langle \mathrm{grad}_c g(z), \zeta-z \rangle] D_z^c D_z^c \left[ c(\zeta,z,t) \mathcal{K}_{a,b}^S(\zeta,z) \right] dt \in \Lambda^{\frac{\alpha}{2}}(\mathcal{U}). \qquad (77)$$

For fixed $w \in \mathcal{U}$ and $z \in \mathcal{U}$ such that $|z-w| \leq \delta$ we consider then the neighborhoods $W(w,\epsilon,\sqrt{\delta})$ and $W'(w,\epsilon,\sqrt{\delta})$, and the function $\psi_w(\zeta)$ from (60) and represent integral from (77) as

$$\int_{\mathcal{U}(\epsilon)\times[0,1]} [\widetilde{g}(\zeta) - g(z) - \langle \mathrm{grad}_c g(z), \zeta-z \rangle] D_z^c D_z^c \left[ c(\zeta,z,t) \mathcal{K}_{a,b}^S(\zeta,z) \right] dt \qquad (78)$$
$$= \int_{\mathcal{U}(\epsilon)\times[0,1]} [\widetilde{g}(\zeta) - g(z) - \langle \mathrm{grad}_c g(z), \zeta-z \rangle] \psi_w(\zeta) D_z^c D_z^c \left[ c(\zeta,z,t) \mathcal{K}_{a,b}^S(\zeta,z) \right] dt$$
$$+ \int_{\mathcal{U}(\epsilon)\times[0,1]} [\widetilde{g}(\zeta) - g(z) - \langle \mathrm{grad}_c g(z), \zeta-z \rangle] (1-\psi_w(\zeta)) D_z^c D_z^c \left[ c(\zeta,z,t) \mathcal{K}_{a,b}^S(\zeta,z) \right] dt.$$

To estimate integrals in the right hand side of (78) we use estimate (51) and obtain a representation
$$D_z^c D_z^c \left[ c(\zeta, z, t) \mathcal{K}_{a,b}^S(\zeta, z) \right] = \sum_{I,d,j} c_{\{I,d,j\}}(\zeta, z, t) \cdot \mathcal{K}_{d,h}^I(\zeta, z) \tag{79}$$

with $c_{\{I,d,j\}}(\zeta, z, t) = c_{\{I,d,j\}}(\zeta, z, \theta(\zeta), t) \in C^\infty \left( \widetilde{\mathcal{U}}_\zeta \times \widetilde{\mathcal{U}}_z \times \mathbb{S}^{n-1} \times [0,1] \right)$, and indices and multiindices in the right hand side of (79) satisfying
$$\begin{aligned} k\left(\mathcal{K}_{d,h}^I\right) + h - l\left(\mathcal{K}_{d,h}^I\right) &\leq 2n - m, \\ k\left(\mathcal{K}_{d,h}^I\right) + 2h - 2l\left(\mathcal{K}_{d,h}^I\right) &\leq 2n - m + 2. \end{aligned} \tag{80}$$

Then for the first integral in the right hand side of (78) using the estimate
$$|\tilde{g}(\zeta) - g(z) - \langle \mathrm{grad}_c g(z), \zeta - z \rangle| = \|g\|_{\Pi^{1+\alpha}(\mathcal{U})} \cdot \mathcal{O}\left(|\Phi(\zeta, z)|^{\frac{1+\alpha}{2}}\right), \tag{81}$$

representation (79) with kernels satisfying (80) and lemma 3.5 we obtain the estimate
$$\left| \int_{\mathcal{U}(\epsilon)\times[0,1]} [\tilde{g}(\zeta) - g(z) - \langle \mathrm{grad}_c g(z), \zeta - z \rangle] \psi_w(\zeta) D_z^c D_z^c \left[ c(\zeta, z, t) \mathcal{K}_{a,b}^S(\zeta, z) \right] dt \right|$$
$$= \left| \sum_{I,d,j} \int_{\mathcal{U}(\epsilon)\times[0,1]} [\tilde{g}(\zeta) - g(z) - \langle \mathrm{grad}_c g(z), \zeta - z \rangle] \psi_w(\zeta) \left[ c_{\{I,d,j\}}(\zeta, z, t) \cdot \mathcal{K}_{d,h}^I(\zeta, z) \right] dt \right|$$
$$= \|g\|_{\Pi^{1+\alpha}(\mathcal{U})} \cdot \mathcal{O} \left( \int_{W_w'(\epsilon,\sqrt{\delta})} \frac{\overbrace{\wedge_i d\theta_i(\zeta)}^{m-1} \wedge d\sigma_{2n-m}(\zeta)}{|\zeta - z|^{k(\mathcal{K})} \cdot |\Phi(\zeta, z)|^{h(\mathcal{K})-l(\mathcal{K})-\frac{1+\alpha}{2}}} \right)$$
$$= \|g\|_{\Pi^{1+\alpha}(\mathcal{U})} \cdot \mathcal{O}\left( \mathcal{I}_1 \left\{\alpha, k(\mathcal{K}), h(\mathcal{K}) - l(\mathcal{K}) - 1/2\right\}\left(\epsilon, \sqrt{\delta}\right) \right)$$
$$= \|g\|_{\Pi^{1+\alpha}(\mathcal{U})} \cdot \mathcal{O}\left(\delta^{\frac{\alpha}{2}}\right).$$

For the second integral in the right hand side of (78) we consider the following representation
$$\int_{\mathcal{U}(\epsilon)\times[0,1]} [\tilde{g}(\zeta) - g(z) - \langle \mathrm{grad}_c g(z), \zeta - z \rangle] (1 - \psi_w(\zeta)) D_z^c D_z^c \left[ c(\zeta, z, t) \mathcal{K}_{a,b}^S(\zeta, z) \right] dt \tag{82}$$
$$- \int_{\mathcal{U}(\epsilon)\times[0,1]} [\tilde{g}(\zeta) - g(w) - \langle \mathrm{grad}_c g(w), \zeta - w \rangle] (1 - \psi_w(\zeta)) D_w^c D_w^c \left[ c(\zeta, w, t) \mathcal{K}_{a,b}^S(\zeta, w) \right] dt$$
$$= \int_{\mathcal{U}(\epsilon)\times[0,1]} [\tilde{g}(\zeta) - g(z) - \langle \mathrm{grad}_c g(z), \zeta - z \rangle] (1 - \psi_w(\zeta))$$
$$\times \left( D_z^c D_z^c \left[ c(\zeta, z, t) \mathcal{K}_{a,b}^S(\zeta, z) \right] - D_w^c D_w^c \left[ c(\zeta, w, t) \mathcal{K}_{a,b}^S(\zeta, w) \right] \right) dt$$
$$+ \int_{\mathcal{U}(\epsilon)\times[0,1]} \langle (\mathrm{grad}_c g(w) - \mathrm{grad}_c g(z)), \zeta - w \rangle (1 - \psi_w(\zeta)) D_w^c D_w^c \left[ c(\zeta, w, t) \mathcal{K}_{a,b}^S(\zeta, w) \right] dt$$
$$+ \int_{\mathcal{U}(\epsilon)\times[0,1]} [g(w) - g(z) - \langle \mathrm{grad}_c g(z), w - z \rangle] (1 - \psi_w(\zeta)) D_w^c D_w^c \left[ c(\zeta, w, t) \mathcal{K}_{a,b}^S(\zeta, w) \right] dt.$$

For the first integral in the right hand side of (82) using estimates (44) and (81), inequalities (80) and lemma 3.5 we obtain
$$\left| \int_{\mathcal{U}(\epsilon)\times[0,1]} [\tilde{g}(\zeta) - g(z) - \langle \mathrm{grad}_c g(z), \zeta - z \rangle] (1 - \psi_w(\zeta)) \right.$$
$$\left. \times \left( D_z^c D_z^c \left[ c(\zeta, z, t) \mathcal{K}_{a,b}^S(\zeta, z) \right] - D_w^c D_w^c \left[ c(\zeta, w, t) \mathcal{K}_{a,b}^S(\zeta, w) \right] \right) dt \right|$$

$$= \left| \sum_{I,d,j} \int_{\mathcal{U}(\epsilon) \times [0,1]} [\tilde{g}(\zeta) - g(z) - \langle \mathrm{grad}_c g(z), \zeta - z \rangle] (1 - \psi_w(\zeta)) \right.$$
$$\left. \times \left[ c_{\{I,d,j\}}(\zeta, z, t) \cdot \mathcal{K}_{d,h}^I(\zeta, z) - c_{\{I,d,j\}}(\zeta, w, t) \cdot \mathcal{K}_{d,h}^I(\zeta, w) \right] dt \right|$$
$$= \|g\|_{\Pi^{1+\alpha}(\mathcal{U})} \cdot \delta \cdot \left[ \mathcal{O}\left(\mathcal{I}_2\{\alpha, k(\mathcal{K}), h(\mathcal{K}) - l(\mathcal{K}) - 1/2\}\left(\epsilon, \sqrt{\delta}\right)\right) \right.$$
$$+ \mathcal{O}\left(\mathcal{I}_2\{\alpha, k(\mathcal{K}) + 1, h(\mathcal{K}) - l(\mathcal{K}) - 1/2\}\left(\epsilon, \sqrt{\delta}\right)\right)$$
$$\left. + \mathcal{O}\left(\mathcal{I}_2\{\alpha, k(\mathcal{K}), h(\mathcal{K}) - l(\mathcal{K}) + 1/2\}\left(\epsilon, \sqrt{\delta}\right)\right) \right]$$
$$= \|g\|_{\Pi^{1+\alpha}(\mathcal{U})} \cdot \delta \cdot \mathcal{O}\left(\delta^{\frac{\alpha-2}{2}}\right) = \|g\|_{\Pi^{1+\alpha}(\mathcal{U})} \cdot \mathcal{O}\left(\delta^{\alpha/2}\right).$$

For the second integral in the right hand side of (82) using the estimate
$$|\mathrm{grad}_c g(w) - \mathrm{grad}_c g(z)| = \|g\|_{\Pi^{1+\alpha}(\mathcal{U})} \cdot \mathcal{O}\left(\delta^{\alpha/2}\right)$$
we reduce the problem to estimate
$$\left| \int_{\mathcal{U}(\epsilon) \times [0,1]} (\zeta_i - w_i)(1 - \psi_w(\zeta)) D_w^c D_w^c \left[ c(\zeta, w, t) \mathcal{K}_{a,b}^S(\zeta, w) \right] dt \right| = \mathcal{O}(1),$$
which follows from the arguments presented in the proof of the second part of lemma 3.7.

For the third integral of the right hand side of (82) using estimate
$$|g(w) - g(z) - \langle \mathrm{grad}_c g(z), w - z \rangle| = \|g\|_{\Pi^{1+\alpha}(\mathcal{U})} \cdot \mathcal{O}\left(\delta^{\frac{1+\alpha}{2}}\right),$$
representation (79) and lemma 3.5 we obtain
$$\left| \int_{\mathcal{U}(\epsilon) \times [0,1]} [g(w) - g(z) - \langle \mathrm{grad}_c g(z), w - z \rangle] (1 - \psi_w(\zeta)) D_w^c D_w^c \left[ c(\zeta, w, t) \mathcal{K}_{a,b}^S(\zeta, w) \right] dt \right|$$
$$= \left| \sum_{I,d,j} \int_{\mathcal{U}(\epsilon) \times [0,1]} [g(w) - g(z) - \langle \mathrm{grad}_c g(z), w - z \rangle] (1 - \psi_w(\zeta)) \right.$$
$$\left. \times \left[ c_{\{I,d,j\}}(\zeta, z, t) \cdot \mathcal{K}_{d,h}^I(\zeta, z) \right] dt \right|$$
$$= \|g\|_{\Pi^{1+\alpha}(\mathcal{U})} \cdot \delta^{\frac{1+\alpha}{2}} \cdot \mathcal{O}\left(\mathcal{I}_2\{0, k(\mathcal{K}), h(\mathcal{K}) - l(\mathcal{K})\}\left(\epsilon, \sqrt{\delta}\right)\right)$$
$$= \|g\|_{\Pi^{1+\alpha}(\mathcal{U})} \cdot \mathcal{O}\left(\delta^{\alpha/2}\right).$$

To prove the inclusion $D^c D^c f_\epsilon \in \Gamma^\alpha(\mathcal{U})$ we consider representation (76) and reduce the problem to the proof of inclusion
$$\int_{\mathcal{U}(\epsilon) \times [0,1]} [\tilde{g}(\zeta) - g(z) - \langle \mathrm{grad}_c g(z), \zeta - z \rangle] D_z^c D_z^c \left[ c(\zeta, z, t) \mathcal{K}_{a,b}^S(\zeta, z) \right] dt \in \Gamma^\alpha(\mathcal{U}). \tag{83}$$

For $z, w \in \mathcal{U}$ such that there exists a curve $x : [0,1] \to \mathbf{M}$ satisfying (2) with $x(0) = z, x(1) = w$, we represent integral from (83) as
$$\int_{\mathcal{U}(\epsilon) \times [0,1]} [\tilde{g}(\zeta) - g(z) - \langle \mathrm{grad}_c g(z), \zeta - z \rangle] D_z^c D_z^c \left[ c(\zeta, z, t) \mathcal{K}_{a,b}^S(\zeta, z) \right] dt \tag{84}$$
$$= \int_{\mathcal{U}(\epsilon) \times [0,1]} [\tilde{g}(\zeta) - g(z) - \langle \mathrm{grad}_c g(z), \zeta - z \rangle] \phi_w(\zeta) D_z^c D_z^c \left[ c(\zeta, z, t) \mathcal{K}_{a,b}^S(\zeta, z) \right] dt$$
$$+ \int_{\mathcal{U}(\epsilon) \times [0,1]} [\tilde{g}(\zeta) - g(z) - \langle \mathrm{grad}_c g(z), \zeta - z \rangle] (1 - \phi_w(\zeta)) D_z^c D_z^c \left[ c(\zeta, z, t) \mathcal{K}_{a,b}^S(\zeta, z) \right] dt$$
using the neighborhoods $W(w, \epsilon, \delta)$ and $W'(w, \epsilon, \delta)$ with $|z - w| \leq \delta$ and a function $\phi_w(\zeta)$ from (49).

For the first integral in (84) we use estimate (81), representation (79) with indices satisfying (80) and lemma 3.5 and obtain the following estimate

$$\left| \int_{\mathcal{U}(\epsilon)\times[0,1]} [\tilde g(\zeta) - g(z) - \langle \mathrm{grad}_c g(z), \zeta - z\rangle ] \phi_w(\zeta) D_z^c D_z^c \left[ c(\zeta, z, t) \mathcal{K}_{a,b}^S(\zeta, z)\right] dt \right|$$

$$= \left| \sum_{I,d,j} \int_{\mathcal{U}(\epsilon)\times[0,1]} [\tilde g(\zeta) - g(z) - \langle \mathrm{grad}_c g(z), \zeta - z\rangle] \phi_w(\zeta) \left[ c_{\{I,d,j\}}(\zeta, z, t) \cdot \mathcal{K}_{d,h}^I(\zeta, z)\right] dt \right|$$

$$= \|g\|_{\Pi^{1+\alpha}(\mathcal{U})} \cdot \mathcal{O}\left( \int_{W'_w(\epsilon,\delta)} \frac{\overbrace{\wedge_i d\theta_i(\zeta)}^{m-1} \wedge d\sigma_{2n-m}(\zeta)}{|\zeta - z|^{k(\mathcal{K})} \cdot |\Phi(\zeta, z)|^{h(\mathcal{K}) - l(\mathcal{K}) - \frac{1+\alpha}{2}}} \right)$$

$$= \|g\|_{\Pi^{1+\alpha}(\mathcal{U})} \cdot \mathcal{O}\left( \mathcal{I}_1\{\alpha, k(\mathcal{K}), h(\mathcal{K}) - l(\mathcal{K}) - 1/2\}(\epsilon, \delta)\right) = \|g\|_{\Pi^{1+\alpha}(\mathcal{U})} \cdot \mathcal{O}(\delta^\alpha).$$

For the second integral in the right hand side of (84) we consider the following representation which is an analogue of representation (82) but with function $\phi_w(\zeta)$ from (49):

$$\int_{\mathcal{U}(\epsilon)\times[0,1]} [\tilde g(\zeta) - g(z) - \langle \mathrm{grad}_c g(z), \zeta - z\rangle] (1 - \phi_w(\zeta)) D_z^c D_z^c \left[ c(\zeta, z, t) \mathcal{K}_{a,b}^S(\zeta, z)\right] dt \quad (85)$$

$$- \int_{\mathcal{U}(\epsilon)\times[0,1]} [\tilde g(\zeta) - g(w) - \langle \mathrm{grad}_c g(w), \zeta - w\rangle] (1 - \phi_w(\zeta)) D_w^c D_w^c \left[ c(\zeta, w, t) \mathcal{K}_{a,b}^S(\zeta, w)\right] dt$$

$$= \int_{\mathcal{U}(\epsilon)\times[0,1]} \langle (\mathrm{grad}_c g(w) - \mathrm{grad}_c g(z)), \zeta - w\rangle (1 - \phi_w(\zeta)) D_w^c D_w^c \left[ c(\zeta, w, t) \mathcal{K}_{a,b}^S(\zeta, w)\right] dt$$

$$+ \int_{\mathcal{U}(\epsilon)\times[0,1]} [g(w) - g(z) - \langle \mathrm{grad}_c g(z), w - z\rangle] (1 - \phi_w(\zeta)) D_w^c D_w^c \left[ c(\zeta, w, t) \mathcal{K}_{a,b}^S(\zeta, w)\right] dt$$

$$+ \int_{\mathcal{U}(\epsilon)\times[0,1]} [\tilde g(\zeta) - g(z) - \langle \mathrm{grad}_c g(z), \zeta - z\rangle] (1 - \phi_w(\zeta))$$

$$\times \left( D_z^c D_z^c \left[ c(\zeta, z, t)\mathcal{K}_{a,b}^S(\zeta, z)\right] - D_w^c D_w^c \left[ c(\zeta, w, t)\mathcal{K}_{a,b}^S(\zeta, w)\right] \right) dt.$$

For the first integral in the right hand side of (85) using the estimate

$$|\mathrm{grad}_c g(w) - \mathrm{grad}_c g(z)| = \|g\|_{\Pi^{1+\alpha}(\mathcal{U})} \cdot \mathcal{O}(\delta^\alpha)$$

and the arguments from the proof of the second part of lemma 3.7 for estimate

$$\left| \int_{\mathcal{U}(\epsilon)\times[0,1]} (\zeta_i - w_i)(1 - \phi_w(\zeta)) D_w^c D_w^c \left[ c(\zeta, w, t) \mathcal{K}_{a,b}^S(\zeta, w)\right] dt \right| = \mathcal{O}(1)$$

we obtain an estimate $\|g\|_{\Pi^{1+\alpha}(\mathcal{U})} \cdot \mathcal{O}(\delta^\alpha)$.

For the second integral of the right hand side of (85) using the estimate

$$|g(w) - g(z) - \langle \mathrm{grad}_c g(z), w - z\rangle| = \|g\|_{\Pi^{1+\alpha}(\mathcal{U})} \cdot \mathcal{O}\left(\delta^{1+\alpha}\right),$$

representation (79) with indices satisfying (80) and lemma 3.5 we obtain an estimate $\|g\|_{\Pi^{1+\alpha}(\mathcal{U})} \cdot \mathcal{O}(\delta^\alpha)$.

To estimate the third integral in the right hand side of (85) we have to modify representation (79) in such a way that the kernels of a new representation will satisfy conditions

$$\begin{aligned} k\left(\mathcal{K}_{d,h}^I\right) + h - l\left(\mathcal{K}_{d,h}^I\right) &\leq 2n - m - 1, \\ k\left(\mathcal{K}_{d,h}^I\right) + 2h - 2l\left(\mathcal{K}_{d,h}^I\right) &\leq 2n - m + 2. \end{aligned} \quad (86)$$

To make such a modification we use the formula

$$D^c_{\bar{z}} D^c_z \left[ c(\zeta, z, t) \cdot \frac{(\zeta - z)^{S_2} (\bar{\zeta} - \bar{z})^{S_3}}{|\zeta - z|^a \cdot \Phi(\zeta, z)^b} \right] \tag{87}$$

$$= D^c_{\bar{z}} \left[ D^c_z \left[ \frac{1}{\Phi(\zeta, z)^b} \right] \cdot \frac{c(\zeta, z, t)(\zeta - z)^{S_2} (\bar{\zeta} - \bar{z})^{S_3}}{|\zeta - z|^a} \right]$$

$$+ D^c_{\bar{z}} \left[ \left( D^c_z + D^c_\zeta \right) \left[ \frac{c(\zeta, z, t)(\zeta - z)^{S_2} (\bar{\zeta} - \bar{z})^{S_3}}{|\zeta - z|^a} \right] \cdot \frac{1}{\Phi(\zeta, z)^b} \right]$$

$$- D^c_\zeta \left[ \frac{c(\zeta, z, t)(\zeta - z)^{S_2} (\bar{\zeta} - \bar{z})^{S_3}}{|\zeta - z|^a} \right] \cdot D^c_{\bar{z}} \left[ \frac{1}{\Phi(\zeta, z)^b} \right]$$

$$- D^c_\zeta D^c_{\bar{z}} \left[ \frac{c(\zeta, z, t)(\zeta - z)^{S_2} (\bar{\zeta} - \bar{z})^{S_3}}{|\zeta - z|^a} \right] \cdot \frac{1}{\Phi(\zeta, z)^b}.$$

Applying estimates (25) and (51) and using inequalities (42) for the kernels $\mathcal{K}^S_{a,b}(\zeta, z)$ we conclude that the first three terms of the right hand side of (87) can be represented as sums of terms of the form $c(\zeta, z, t) \cdot \mathcal{K}^I_{d,h}(\zeta, z)$ with indices satisfying (86).

For such kernels using estimates (55) and (81), and lemma 3.5 we obtain in the last integral of the right hand side of (85)

$$\left| \int_{\mathcal{U}(\epsilon) \times [0,1]} [\tilde{g}(\zeta) - g(z) - \langle \mathrm{grad}_c g(z), \zeta - z \rangle] (1 - \phi_w(\zeta)) \right.$$

$$\left. \times \left[ c(\zeta, z, t) \cdot \mathcal{K}^I_{d,h}(\zeta, z) - c(\zeta, w, t) \cdot \mathcal{K}^I_{d,h}(\zeta, w) \right] dt \right|$$

$$= \|g\|_{\Pi^{1+\alpha}(\mathcal{U})} \cdot \left[ \delta \cdot \mathcal{O} \left( \mathcal{I}_2 \{\alpha, k(\mathcal{K}), h(\mathcal{K}) - l(\mathcal{K}) - 1/2\} (\epsilon, \delta) \right) \right.$$

$$+ \delta \cdot \mathcal{O} \left( \mathcal{I}_2 \{\alpha, k(\mathcal{K}) + 1, h(\mathcal{K}) - l(\mathcal{K}) - 1/2\} (\epsilon, \delta) \right)$$

$$+ \delta \cdot \mathcal{O} \left( \mathcal{I}_2 \{\alpha, k(\mathcal{K}) - 1, h(\mathcal{K}) - l(\mathcal{K}) + 1/2\} (\epsilon, \delta) \right)$$

$$\left. + \delta^2 \cdot \mathcal{O} \left( \mathcal{I}_2 \{\alpha, k(\mathcal{K}), h(\mathcal{K}) - l(\mathcal{K}) + 1/2\} (\epsilon, \delta) \right) \right]$$

$$= \|g\|_{\Pi^{1+\alpha}(\mathcal{U})} \cdot \left[ \delta \cdot \mathcal{O} \left( \delta^{\alpha - 1} \right) + \delta^2 \cdot \mathcal{O} \left( \delta^{\alpha - 2} \right) \right] = \|g\|_{\Pi^{1+\alpha}(\mathcal{U})} \cdot \mathcal{O} \left( \delta^\alpha \right).$$

For the last kernel from the right hand side of (87) integrating by parts we obtain

$$\int_{\mathcal{U}(\epsilon) \times [0,1]} [\tilde{g}(\zeta) - g(z) - \langle \mathrm{grad}_c g(z), \zeta - z \rangle] (1 - \phi_w(\zeta)) \tag{88}$$

$$\times \left[ D^c_\zeta D^c_{\bar{z}} \left[ \frac{c(\zeta, z, t)(\zeta - z)^{S_2} (\bar{\zeta} - \bar{z})^{S_3}}{|\zeta - z|^a} \right] \cdot \frac{1}{\Phi(\zeta, z)^b} \right.$$

$$\left. - D^c_\zeta D^c_{\bar{w}} \left[ \frac{c(\zeta, w, t)(\zeta - w)^{S_2} (\bar{\zeta} - \bar{w})^{S_3}}{|\zeta - w|^a} \right] \cdot \frac{1}{\Phi(\zeta, w)^b} \right] dt$$

$$= \int_{\mathcal{U}(\epsilon) \times [0,1]} [\tilde{g}(\zeta) - g(z) - \langle \mathrm{grad}_c g(z), \zeta - z \rangle] (1 - \phi_w(\zeta))$$

$$\times \left[ D^c_{\bar{z}} \left[ \frac{c(\zeta, z, t)(\zeta - z)^{S_2} (\bar{\zeta} - \bar{z})^{S_3}}{|\zeta - z|^a} \right] \cdot D^c_\zeta \left[ \frac{1}{\Phi(\zeta, z)^b} \right] \right.$$

$$\left. - D^c_{\bar{w}} \left[ \frac{c(\zeta, w, t)(\zeta - w)^{S_2} (\bar{\zeta} - \bar{w})^{S_3}}{|\zeta - w|^a} \right] \cdot D^c_\zeta \left[ \frac{1}{\Phi(\zeta, w)^b} \right] \right] dt$$

$$+ \int_{\mathcal{U}(\epsilon) \times [0,1]} \left[ D^c_\zeta \tilde{g}(\zeta) - \left\langle \mathrm{grad}_c g(z), D^c_\zeta \zeta \right\rangle \right] (1 - \phi_w(\zeta))$$

$$\times \left[ D_z^c \left[ \frac{c(\zeta,z,t)(\zeta-z)^{S_2}(\bar\zeta-\bar z)^{S_3}}{|\zeta-z|^a} \right] \cdot \frac{1}{\Phi(\zeta,z)^b} \right.$$
$$\left. - D_w^c \left[ \frac{c(\zeta,w,t)(\zeta-w)^{S_2}(\bar\zeta-\bar w)^{S_3}}{|\zeta-w|^a} \right] \cdot \frac{1}{\Phi(\zeta,w)^b} \right] dt$$
$$+ \int_{\mathcal{U}(\epsilon)\times[0,1]} [\widetilde{g}(\zeta) - g(z) - \langle \mathrm{grad}_c g(z), \zeta-z \rangle] \cdot D_\zeta^c(1-\phi_w(\zeta))$$
$$\times \left[ D_z^c \left[ \frac{c(\zeta,z,t)(\zeta-z)^{S_2}(\bar\zeta-\bar z)^{S_3}}{|\zeta-z|^a} \right] \cdot \frac{1}{\Phi(\zeta,z)^b} \right.$$
$$\left. - D_w^c \left[ \frac{c(\zeta,w,t)(\zeta-w)^{S_2}(\bar\zeta-\bar w)^{S_3}}{|\zeta-w|^a} \right] \cdot \frac{1}{\Phi(\zeta,w)^b} \right] dt.$$

Kernels of the first term of the right hand side of (88) satisfy inequalities (86) and therefore admit an estimate $\|g\|_{\Pi^{1+\alpha}(\mathcal{U})} \cdot \mathcal{O}(\delta^\alpha)$.

For the second term of the right hand side of (88) using the estimate
$$\left| D_\zeta^c \widetilde{g}(\zeta) - \left\langle \mathrm{grad}_c g(z), D_\zeta^c \zeta \right\rangle \right| = \|g\|_{\Pi^{1+\alpha}(\mathcal{U})} \cdot \mathcal{O}\left(|\Phi(\zeta,z)|^{\frac{\alpha}{2}}\right)$$
and estimate (55) we obtain
$$\left| \int_{\mathcal{U}(\epsilon)\times[0,1]} \left[ D_\zeta^c \widetilde{g}(\zeta) - \left\langle \mathrm{grad}_c g(z), D_\zeta^c \zeta \right\rangle \right] (1-\phi_w(\zeta)) \right.$$
$$\times \left[ D_z^c \left[ \frac{c(\zeta,z,t)(\zeta-z)^{S_2}(\bar\zeta-\bar z)^{S_3}}{|\zeta-z|^a} \right] \cdot \frac{1}{\Phi(\zeta,z)^b} \right.$$
$$\left.\left. - D_w^c \left[ \frac{c(\zeta,w,t)(\zeta-w)^{S_2}(\bar\zeta-\bar w)^{S_3}}{|\zeta-w|^a} \right] \cdot \frac{1}{\Phi(\zeta,w)^b} \right] dt \right|$$
$$= \|g\|_{\Pi^{1+\alpha}(\mathcal{U})} \cdot [\delta \cdot \mathcal{O}(\mathcal{I}_2\{\alpha, k(\mathcal{K}), h(\mathcal{K}) - l(\mathcal{K})\}(\epsilon,\delta))$$
$$+ \delta \cdot \mathcal{O}(\mathcal{I}_2\{\alpha, k(\mathcal{K})+2, h(\mathcal{K}) - l(\mathcal{K})\}(\epsilon,\delta))$$
$$+ \delta \cdot \mathcal{O}(\mathcal{I}_2\{\alpha, k(\mathcal{K}), h(\mathcal{K}) - l(\mathcal{K})+1\}(\epsilon,\delta))$$
$$+ \delta^2 \cdot \mathcal{O}(\mathcal{I}_2\{\alpha, k(\mathcal{K})+1, h(\mathcal{K}) - l(\mathcal{K})+1\}(\epsilon,\delta))]$$
$$= \|g\|_{\Pi^{1+\alpha}(\mathcal{U})} \cdot \left[ \delta \cdot \mathcal{O}\left(\delta^{\alpha-1}\right) + \delta^2 \cdot \mathcal{O}\left(\delta^{\alpha-2}\right) \right] = \|g\|_{\Pi^{1+\alpha}(\mathcal{U})} \cdot \mathcal{O}(\delta^\alpha).$$

For the third term of the right hand side of (88) using estimates (49), (55), (81) and lemma 3.5 we obtain
$$\left| \int_{\mathcal{U}(\epsilon)\times[0,1]} [\widetilde{g}(\zeta) - g(z) - \langle \mathrm{grad}_c g(z), \zeta-z \rangle] \cdot D_\zeta^c(1-\phi_w(\zeta)) \right.$$
$$\times \left[ D_z^c \left[ \frac{c(\zeta,z,t)(\zeta-z)^{S_2}(\bar\zeta-\bar z)^{S_3}}{|\zeta-z|^a} \right] \cdot \frac{1}{\Phi(\zeta,z)^b} \right.$$
$$\left.\left. - D_w^c \left[ \frac{c(\zeta,w,t)(\zeta-w)^{S_2}(\bar\zeta-\bar w)^{S_3}}{|\zeta-w|^a} \right] \cdot \frac{1}{\Phi(\zeta,w)^b} \right] dt \right|$$
$$= \|g\|_{\Pi^{1+\alpha}(\mathcal{U})} \cdot [\mathcal{O}(\mathcal{I}_1\{\alpha, k(\mathcal{K}), h(\mathcal{K}) - l(\mathcal{K}) - 1/2\}(\epsilon,\delta))$$
$$+ \mathcal{O}(\mathcal{I}_1\{\alpha, k(\mathcal{K})+2, h(\mathcal{K}) - l(\mathcal{K}) - 1/2\}(\epsilon,\delta))$$
$$+ \mathcal{O}(\mathcal{I}_1\{\alpha, k(\mathcal{K}), h(\mathcal{K}) - l(\mathcal{K}) + 1/2\}(\epsilon,\delta))$$
$$+ \delta \cdot \mathcal{O}(\mathcal{I}_2\{\alpha, k(\mathcal{K})+1, h(\mathcal{K}) - l(\mathcal{K}) + 1/2\}(\epsilon,\delta))]$$
$$= \|g\|_{\Pi^{1+\alpha}(\mathcal{U})} \cdot \left[ \mathcal{O}(\delta^\alpha) + \delta \cdot \mathcal{O}\left(\delta^{\alpha-1}\right) \right] = \|g\|_{\Pi^{1+\alpha}(\mathcal{U})} \cdot \mathcal{O}(\delta^\alpha).$$

This concludes the proof of inclusion $D^c D^c f_\epsilon \in \Gamma^\alpha(\mathcal{U})$. The proof of inclusion $D f_\epsilon \in \Gamma^\alpha(\mathcal{U})$ is completely analogous because kernels $D\left[c(\zeta, z, t) \mathcal{K}_{a,b}^S(\zeta, z)\right]$ admit representation (79) with indices satisfying inequalities (86). $\square$

In order to prove applicability of lemmas 3.7 and 3.8 to the kernels obtained from $\lambda_{r-1}^{i,J}$ and $\gamma_{r-1}^{i,J}$ after applications of lemmas 3.2, 3.3 we have to prove relations (42) for these kernels. But according to lemmas 3.2, 3.3 expressions in the left hand sides of these relations don't increase under transformations from these lemmas. Therefore it suffices to prove relations (42) for the original kernels $\mathcal{K}_{d,h}^I(\zeta, z)$ satisfying conditions (34) and (35).

Second condition from (42) is always satisfied for the indices satisfying (34) as can be seen from the inequality
$$k(\mathcal{K}) + 2h(\mathcal{K}) - 2l(\mathcal{K}) \leq 2n - m - |J_6| \leq 2n - m, \tag{89}$$
where we used lemma 3.4 and relations
$$|J_1| + |J_4| + r + m - n \geq 0,$$
$$\sum_{i=1}^4 |J_i| = n - r - 1,$$
$$|J_2| + |J_3| \leq m - 1$$
for the multiindices of $\lambda_{r-1}^{i,J}$.

The same arguments show that condition (89) is also satisfied for the indices defined by (35).

First condition from (42) is not satisfied for all kernels $\mathcal{K}_{d,h}^I(\zeta, z)$. But in the lemma below we show that if this condition is not satisfied then the corresponding term of the integral formula for $R_r(\epsilon)$ does not survive under the limit when $\epsilon \to 0$.

**Lemma 3.9.** *If $k(\mathcal{K}), h(\mathcal{K}), l(\mathcal{K}) \in \mathbb{Z}$ and*
$$k(\mathcal{K}) + h(\mathcal{K}) - l(\mathcal{K}) \geq 2n - m - 1$$
*then*
$$\left\| \int_{\mathcal{U}(\epsilon) \times [0,1]} \widetilde{g}(\zeta) c(\zeta, z, t) \mathcal{K}_{d,h}^I(\zeta, z) dt \right\|_{L^\infty(\mathbf{M})} = \mathcal{O}(\sqrt{\epsilon} \cdot \log \epsilon) \cdot \|g\|_{L^\infty(\mathbf{M})}.$$

**Proof.** We use inequality
$$2n - m + l(\mathcal{K}) - k(\mathcal{K}) - h(\mathcal{K}) \geq n - |J_1| - |J_5| + |J_6| - 1 \geq 1, \tag{90}$$
which is a corollary of definitions of $k(\mathcal{K})$, $h(\mathcal{K})$ and $l(\mathcal{K})$ and equality
$$\sum_{i=1}^4 |J_i| = n - r - 1.$$

From the condition of the lemma and inequality (90) we obtain
$$k(\mathcal{K}) + h(\mathcal{K}) - l(\mathcal{K}) = 2n - m - 1$$
and
$$n - |J_1| - |J_5| - 1 = 1,$$
which leads to
$$|J_1| = n - r - 1, \quad |J_3| = 0, \quad |J_4| = 0, \quad |J_5| = r - 1,$$
and
$$l(\mathcal{K}) \geq |J_1| + |J_4| + r + m - n \geq m - 1 \geq 1.$$

Using lemma 3.5 to estimate the integral in lemma we obtain

$$\left| \int_{\mathcal{U}(\epsilon) \times [0,1]} \widetilde{g}(\zeta) c(\zeta, z, t) \mathcal{K}_{d,h}^I(\zeta, z) dt \right|$$

$$= \|g\|_{L^\infty(\mathbf{M})} \cdot \epsilon^{l(\mathcal{K})} \cdot \mathcal{O}\left(\mathcal{I}_1\{0, k(\mathcal{K}), h(\mathcal{K}), 0\}(\epsilon, 1)\right)$$

$$= \|g\|_{L^\infty(\mathbf{M})} \cdot \begin{cases} \epsilon^{l(\mathcal{K})} \cdot \mathcal{O}\left(\epsilon^{2n-m-k(\mathcal{K})-h(\mathcal{K})} \cdot (\log \epsilon)^2\right) & \text{if } k(\mathcal{K}) \geq 2n - m - 1, \\ \epsilon^{l(\mathcal{K})} \cdot \mathcal{O}\left(\epsilon^{(2n-m-k(\mathcal{K})-2h(\mathcal{K})+1)/2} \cdot \log \epsilon\right) & \text{if } k(\mathcal{K}) \leq 2n - m - 2, \end{cases}$$

where in the first subcase of the above we have the necessary estimate because of inequality (90) and in the second subcase we have the necessary estimate from inequality (89). $\square$

This completes the proof of proposition 3.1.

## 4. Compactness of $\mathbf{H}_r$.

From the definition of operator $\mathbf{H}_r$ we conclude that in order to prove its compactness it suffices to prove compactness of each of the terms below

$$\bar{\partial}_{\mathbf{M}} \vartheta'_\iota(z) \wedge R_r^\iota(\vartheta_\iota g)(z), \quad \vartheta'_\iota(z) \cdot R_{r+1}^\iota(\bar{\partial}_{\mathbf{M}} \vartheta_\iota \wedge g)(z) \quad \text{and} \quad \vartheta'_\iota(z) \cdot H_r^\iota(\vartheta_\iota g)(z).$$

Compactness of the first two of these terms follows from the boundedness of operators $R_r$ proved in proposition 3.1 and compactness of the embedding

$$\Lambda^a(\mathcal{U}) \to \Lambda^b(\mathcal{U})$$

for $a > b$ [Ad]. The proposition below takes care of the third term of $\mathbf{H}_r$.

**Proposition 4.1.** *Let $r < q$. Then*

$$H_r(g)(z) = 0. \tag{91}$$

**Proof.** Using approximation of $H_r$ by the operators

$$H_r(\epsilon)(g)(z) = (-1)^r \frac{(n-1)!}{(2\pi i)^n} \cdot \text{pr}_{\mathbf{M}} \circ \int_{\mathbf{M}_\epsilon} \vartheta(\zeta) \widetilde{g}(\zeta) \wedge \omega'_r\left(\frac{P(\zeta, z)}{\Phi(\zeta, z)}\right) \wedge \omega(\zeta)$$

we conclude that it suffices to prove equality

$$\omega'_r\left(\frac{P(\zeta, z)}{\Phi(\zeta, z)}\right) \wedge \omega(\zeta) = 0$$

for $r < q$.

This kernel with the use of (12) may be represented on $\widetilde{\mathcal{U}} \times \mathcal{U}$ as

$$\omega'_r\left(\frac{P(\zeta, z)}{\Phi(\zeta, z)}\right) \wedge \omega(\zeta) \bigg|_{\widetilde{\mathcal{U}} \times \mathcal{U}} \tag{92}$$

$$= \sum_{i,J} a_{(i,J)}(\zeta, z) \wedge \widetilde{\phi}_r^{i,J}(\zeta, z) + \sum_{i,J} b_{(i,J)}(\zeta, z) \wedge \widetilde{\psi}_r^{i,J}(\zeta, z),$$

where $i$ is an index, $J = \cup_{i=1}^6 J_i$ is a multiindex such that $i \notin J$, $a_{(i,J)}(\zeta, z)$ and $b_{(i,J)}(\zeta, z)$ are smooth functions of $z$, $\zeta$ and $\theta(\zeta)$, and $\widetilde{\phi}_r^{i,J}(\zeta, z)$ and $\widetilde{\psi}_r^{i,J}(\zeta, z)$ are defined as follows:

$$\widetilde{\phi}_r^{i,J}(\zeta, z) = \frac{1}{\Phi(\zeta, z)^n} \times \text{Det}\left[Q^{(i)}, \overbrace{\bar{A} \cdot \bar{\partial}_\zeta a}^{j \in J_1}, \overbrace{a \cdot \mu_\nu}^{j \in J_2}, \overbrace{a \cdot \mu_\tau}^{j \in J_3}, \overbrace{\bar{A} \cdot \bar{\partial}_z a}^{j \in J_4}, \overbrace{a \cdot \bar{\partial}_z \bar{A}}^{j \in J_5}, \overbrace{\bar{\partial}_z Q}^{j \in J_6}\right] \wedge \omega(\zeta) \tag{93}$$

and

$$\widetilde{\psi}_r^{i,J}(\zeta,z) = \frac{1}{\Phi(\zeta,z)^n} \times \mathrm{Det}\left[a_i\bar{A}_i, \overbrace{\bar{A}\cdot\bar{\partial}_\zeta a}^{j\in J_1}, \overbrace{a\cdot\mu_\nu}^{j\in J_2}, \overbrace{a\cdot\mu_\tau}^{j\in J_3}, \overbrace{\bar{A}\cdot\bar{\partial}_z a}^{j\in J_4}, \overbrace{a\cdot\bar{\partial}_z\bar{A}}^{j\in J_5}, \overbrace{\bar{\partial}_z Q}^{j\in J_6}\right]\wedge\omega(\zeta). \quad (94)$$

Multiindices of $\widetilde{\phi}_r^{i,J}$ and $\widetilde{\psi}_r^{i,J}$ satisfy the following conditions

$$\begin{aligned}|J_1|+|J_2|+|J_3| &= n-r-1,\\ |J_1|+|J_2| &\leq m-1,\end{aligned} \quad (95)$$

therefore, if $r < q$ then

$$|J_3| = n-r-1-|J_1|-|J_2| \geq n-r-m > n-q-m,$$

which is impossible. $\square$

DEPARTMENT OF MATHEMATICS, UNIVERSITY OF WYOMING, LARAMIE, WY 82071
*E-mail address*: `polyakov@uwyo.edu`